# Convergence Rates Analysis of The Quadratic Penalty Method and Its Applications to Decentralized Distributed Optimization

Huan Li*    Cong Fang†    Zhouchen Lin‡

November 29, 2017

## Abstract

In this paper, we study a variant of the quadratic penalty method for linearly constrained convex problems, which has already been widely used but actually lacks theoretical justification. Namely, the penalty parameter steadily increases and the penalized objective function is minimized inexactly rather than exactly, e.g., with only one step of the proximal gradient descent. For such a variant of the quadratic penalty method, we give counterexamples to show that it may not give a solution to the original constrained problem. By choosing special penalty parameters, we ensure the convergence and further establish the convergence rates of $O\left(\frac{1}{\sqrt{K}}\right)$ for the generally convex problems and $O\left(\frac{1}{K}\right)$ for strongly convex ones, where $K$ is the number of iterations. Furthermore, by adopting Nesterov's extrapolation we show that the convergence rates can be improved to $O\left(\frac{1}{K}\right)$ for the generally convex problems and $O\left(\frac{1}{K^2}\right)$ for strongly convex ones.

When applied to the decentralized distributed optimization, the penalty methods studied in this paper become the widely used distributed gradient method and the fast distributed gradient method. However, due to the totally different analysis framework, we can improve their $O\left(\frac{\log K}{\sqrt{K}}\right)$ and $O\left(\frac{\log K}{K}\right)$ convergence rates to $O\left(\frac{1}{\sqrt{K}}\right)$ and $O\left(\frac{1}{K}\right)$ with fewer assumptions on the network topology for general convex problems. Using our analysis framework, we also extend the fast distributed gradient method to a communication efficient version, i.e., finding an $\varepsilon$ solution in $O\left(\frac{1}{\varepsilon}\right)$ communications and $O\left(\frac{1}{\varepsilon^{2+\delta}}\right)$ computations for the non-smooth problems, where $\delta$ is a small constant.

## 1 Introduction

The problem of interest in this paper is the linearly constrained convex problem:

$$\min_{\mathbf{x}\in\mathbb{R}^d} h(\mathbf{x}) + f(\mathbf{x}), \quad s.t. \quad \mathbf{Ax} = \mathbf{b}, \tag{1}$$

where $f$ is a convex function with Lipschitz continuous gradient: $\|\nabla f(\mathbf{x}) - \nabla f(\mathbf{y})\| \leq L\|\mathbf{x} - \mathbf{y}\|, \forall \mathbf{x}, \mathbf{y}$. $h$ is convex and can be non-smooth. In the extreme case, $f$ can vanish. For brevity, we denote $F(\mathbf{x}) = h(\mathbf{x}) + f(\mathbf{x})$. We assume that there exists an optimal solution $\mathbf{x}^*$ of (1) and there exists $\lambda^*$ such that $(\mathbf{x}^*, \lambda^*)$ is a saddle point of problem $\min_\mathbf{x} \max_\lambda F(\mathbf{x}) + \langle \lambda, \mathbf{Ax} - \mathbf{b} \rangle$. We focus on $\varepsilon$-approximate minimizer which satisfies $|F(\mathbf{x}_K) - F(\mathbf{x}^*)| \leq O(\varepsilon)$ and $\|\mathbf{Ax}_K - \mathbf{b}\| \leq O(\varepsilon)$, where $K$ is the number of iterations, and $\|\cdot\|$ is the $\ell_2$ norm for a vector and the Frobenius norm for a matrix.

*Peking University. Email: lihuanss@pdu.edu.cn
†Peking University. Email: fangcong@pdu.edu.cn
‡Peking University. Email: zlin@pku.edu.cn

The quadratic penalty method is a natural yet relatively old way to solve (1). It considers the following unconstrained problem:

$$\min_{\mathbf{x} \in \mathbb{R}^d} G(\mathbf{x}) = h(\mathbf{x}) + f(\mathbf{x}) + \frac{\beta}{2} \|\mathbf{A}\mathbf{x} - \mathbf{b}\|^2 \tag{2}$$

instead. It is known that the minimizer of (2) is a near optimal solution to (1) if the penalty $\beta$ is large enough, and when $\beta \to +\infty$ the solution to (2) coincides with that of (1) [1] [1].

However, directly solving (2) with a large pre-defined penalty will make the solution to (2) very slow [1]. So in practice, to address this issue the continuation technique is used [1, 3, 4, 5, 6, 7, 8, 9]. Namely, solving a sequence of problem (2) with the penalty parameter $\beta$ gradually increasing from a relatively small value. Though continuation is widely adopted in the penalty method [1, 3, 4, 5, 6, 7, 8, 9], its practical implementation is actually different from its theoretical analysis. Namely, most of the existing theoretical analysis requires solving problem (2) *exactly* for each value of $\beta$ [1, 3, 4]. With such an assumption, the convergence of the solution to (2) to that of (1) can be established, and some convergence rates are also provided [1, 3, 4]. However, in reality people do not solve problem (2) exactly since solving problem (2) with a large $\beta$ exactly is time consuming due to its ill-conditioning. Instead, they only solve problem (2) inexactly, e.g., apply the proximal gradient descent to problem (2) *only once* and then increase $\beta$ immediately. Thus recent work relaxes exact solution to problem (2) to approximate solution [8, 9]. However, the approximation error needs to be carefully controlled. This is still inconsistent with the real practice in many engineering applications, where people tend to use simple stopping criteria when solving problem (2), e.g., iterate for a fixed number of iterations (typically only once, as far as we know from the machine learning community). With such a variant of the quadratic penalty method, whether the iterates converge to the solution to (1) is unsure (please see Section 3 for counterexamples), not to say the convergence rate.

Decentralized distributed computation can be widely found in scientific and engineering areas including automatic control, signal processing and machine learning. For specific applications, please see the references in [10]. Existing decentralized algorithms include the distributed subgradient (gradient) method [11, 12, 13, 14, 15, 16, 17], the fast gradient method [18], the dual averaging method [19], the dual decomposition [20], ADMM [21, 22, 23] and EXTRA [10, 24]. For the general convex problem, [12] shows the $O\left(\frac{\log K}{\sqrt{K}}\right)$ convergence rate for the distributed gradient method and [19] established the $O\left(\frac{\log K}{\sqrt{K}}\right)$ convergence rate for the dual averaging method, which has a special case of the subgradient method and gradient method. [18] proved the convergence rate of $O\left(\frac{\log K}{K}\right)$ for the fast gradient method. Dual decomposition [20], ADMM [23] and EXTRA [10] have the convergence rate of $O\left(\frac{1}{K}\right)$. Recently, [25] proposed a communication efficient distributed algorithm based on the primal-dual method, which can find an $\varepsilon$ solution in $O\left(\frac{1}{\varepsilon}\right)$ communications and $O\left(\frac{1}{\varepsilon^2}\right)$ computations. In comparison, the distributed subgradient method needs $O\left(\frac{\log 1/\varepsilon}{\varepsilon^2}\right)$ communications and computations.

The decentralized distributed optimization can be formulated as a constrained problem [20, 21, 22, 23]. We can see that the proved convergence rates of the distributed gradient method [12, 19] and the fast gradient method [18] are generally weaker than the convergence rates we often see in constrained optimization by a factor of $\log K$, such as the non-ergodic $O\left(\frac{1}{\sqrt{K}}\right)$ convergence rate of ADMM [26] and ergodic $O\left(\frac{1}{K}\right)$ convergence rate for ADMM, the Augmented Lagrangian method and the primal-dual method [27, 28]. An open problem is that can we improve the convergence rates of such widely used methods in distributed community from $O\left(\frac{\log K}{\sqrt{K}}\right)$ and $O\left(\frac{\log K}{K}\right)$ to $O\left(\frac{1}{\sqrt{K}}\right)$ and $O\left(\frac{1}{K}\right)$?

---

[1] For some particular problems, e.g., sparsity penalized problems [2], $\beta \to +\infty$ is not a necessary condition for making the solutions to (2) and (1) identical. This situation is not of interest of our paper.

## 1.1 Contributions

In this paper, we aim at filling in the gap between the theories and the practice of the quadratic penalty method for problem (1). We assume that problem (2) is solved by the Proximal Gradient (PG) method or the Accelerated Proximal Gradient (APG) method [29], two simple first order algorithms widely adopted by the machine learning community. Then in practice people will perform the proximal gradient descent only once and then increase $\beta$. With carefully chosen penalties, for PG we establish the convergence rates of $O\left(\frac{1}{\sqrt{K}}\right)$ and $O\left(\frac{1}{K}\right)$ when $f$ is generally convex and strongly convex, respectively; for APG we improve the convergence rates to $O\left(\frac{1}{K}\right)$ and $O\left(\frac{1}{K^2}\right)$. We further give counterexamples to show that the heuristic methods currently widely used in practice may not converge to the solution to the original constrained problem.

When applied to the decentralized distributed optimization, we found that the PG based penalty method and the APG based penalty method become the distributed gradient method in [11, 19, 13] and the fast distributed gradient method in [18], respectively. Due to the totally different analysis framework, we can improve the $O\left(\frac{\log K}{\sqrt{K}}\right)$ and $O\left(\frac{\log K}{K}\right)$ convergence rates in these literatures to $O\left(\frac{1}{\sqrt{K}}\right)$ and $O\left(\frac{1}{K}\right)$ under fewer assumptions on the network topology. Using our analysis framework, we can also easily extend the method in [18] to a communication efficient version, i.e., the algorithm finds an $\varepsilon$ solution in $O\left(\frac{1}{\varepsilon}\right)$ communications and $O\left(\frac{1}{\varepsilon^{2+\delta}}\right)$ computations for the non-smooth problems, where $\delta$ is a small constant.

## 2 The Quadratic Penalty Method with Continuation and Inexact Update

In this section, we study the quadratic penalty method with continuation and inexact update. Section 2.1 derives the convergence rate when PG is used to minimize (2) while Section 2.2 studies the case of APG.

### 2.1 PG as the Solver

In this subsection, we use PG to solve problem (2) inexactly and then increase the penalty. Actually, in each iteration we solve (2) with only one iteration, rather than solving (2) exactly or approximately up to certain precisions as the existing literatures assume. Such a treatment fits for real practice better. Concretely, at each iteration the method consists of the following standard proximal gradient step:

$$\mathbf{x}_{k+1} = \text{prox}_{h,\epsilon_k}^{1/\eta_k}\left(\mathbf{x}_k - \frac{\nabla f(\mathbf{x}_k) + \frac{\beta}{\alpha_k}\mathbf{A}^T(\mathbf{A}\mathbf{x}_k - \mathbf{b})}{\eta_k}\right), \tag{3}$$

which is obtained by the approximate minimization of

$$\begin{aligned}\min_{\mathbf{x}} \quad & h(\mathbf{x}) + f(\mathbf{x}_k) + \langle \nabla f(\mathbf{x}_k), \mathbf{x} - \mathbf{x}_k \rangle + \frac{L}{2}\|\mathbf{x} - \mathbf{x}_k\|^2 \\ & + \frac{\beta}{\alpha_k}\left\langle \mathbf{A}^T(\mathbf{A}\mathbf{x}_k - \mathbf{b}), \mathbf{x} - \mathbf{x}_k \right\rangle + \frac{\beta\|\mathbf{A}^T\mathbf{A}\|_2}{2\alpha_k}\|\mathbf{x} - \mathbf{x}_k\|^2,\end{aligned}$$

where $\eta_k = L + \frac{\beta\|\mathbf{A}^T\mathbf{A}\|_2}{\alpha_k}$ and $\mathbf{x}_{k+1} = \text{prox}_{h,\epsilon_k}^{1/\eta_k}(\mathbf{z})$ means that $h(\mathbf{x}_{k+1}) + \frac{\eta_k}{2}\|\mathbf{x}_{k+1} - \mathbf{z}\|^2 \leq \min_{\mathbf{x}} h(\mathbf{x}) + \frac{\eta_k}{2}\|\mathbf{x} - \mathbf{z}\|^2 + \epsilon_k^2$. We can set $\epsilon_k = 0$ when the proximal mapping $\text{prox}_h^{1/\eta_k}(\cdot)$ has closed form solution.

---

[2]For some cases, i.e., the $l_p$ norm, although its proximal mapping has no closed form solution, it is much more efficient to find its approximate proximal mapping [30] than finding an approximate solution to problem (2).

Comparing with (2), for notational simplicity we have replaced the penalty $\beta$ with $\frac{\beta}{\alpha_k}$, where $\alpha_k \geq 0$ is decreasing, to show that the actual penalty increases. The inexactness comes from two aspects: 1. We linearize $f(\mathbf{x}) + \frac{\beta}{2\alpha_k}\|\mathbf{A}\mathbf{x} - \mathbf{b}\|^2$ in (2). 2. We allow to compute the proximal mapping of $h(\mathbf{x})$ approximately when such proximal mapping has no closed form solution.

We first give a general result in Theorem 1, which considers both the generally convex case and the strongly convex case.

**Theorem 1** *Assume that $f(\mathbf{x})$ is strongly convex with modulus $\mu \geq 0$. Let $\alpha_k \geq 0$ be a decreasing sequence with $\frac{1}{\alpha_{-1}} = 0$. Define $\theta_k$ as $\frac{1-\theta_k}{\alpha_k} = \frac{1}{\alpha_{k-1}}$. Assume that*

$$\frac{(\eta_k - \mu)\theta_k}{2\alpha_k} \leq \frac{\eta_{k-1}\theta_{k-1}}{2\alpha_{k-1}}, \tag{4}$$

$$\sum_{k=0}^{\infty} \frac{\epsilon_k}{\alpha_k} < \infty, \quad \sum_{k=0}^{\infty} \left( \frac{2\sqrt{\epsilon_{k+1}\eta_{k+1}\alpha_k}}{\alpha_{k+1}\sqrt{\eta_k\theta_k}} + \frac{2\sqrt{\epsilon_k}}{\sqrt{\theta_k\alpha_k}} \right) < \infty. \tag{5}$$

*Then for the PG based quadratic penalty method (3), we have*

$$|F(\mathbf{x}_{K+1}) - F(\mathbf{x}^*)| \leq O(\alpha_K), \quad \|\mathbf{A}\mathbf{x}_{K+1} - \mathbf{b}\| \leq O(\alpha_K),$$

*where $\{\mathbf{x}^*\}$ is an optimal solution to (1).*

**Proof 1** *Let $\overline{\boldsymbol{\lambda}}_{k+1} = \frac{\beta}{\alpha_k}(\mathbf{A}\mathbf{x}_k - \mathbf{b})$. Then from the $\epsilon_k$ optimality condition of (3) [31], we have that there exists $\sigma_k$ such that $\|\sigma_k\| \leq \sqrt{\frac{2\epsilon_k}{\eta_k}}$ and*

$$0 \in \nabla f(\mathbf{x}_k) + \partial_{\epsilon_k} h(\mathbf{x}_{k+1}) + \mathbf{A}^T \overline{\boldsymbol{\lambda}}_{k+1} + \eta_k(\mathbf{x}_{k+1} - \mathbf{x}_k + \sigma_k).$$

*So*

$$h(\mathbf{x}) - h(\mathbf{x}_{k+1}) \geq -\epsilon_k - \left\langle \nabla f(\mathbf{x}_k) + \mathbf{A}^T \overline{\boldsymbol{\lambda}}_{k+1} + \eta_k(\mathbf{x}_{k+1} - \mathbf{x}_k + \sigma_k), \mathbf{x} - \mathbf{x}_{k+1} \right\rangle.$$

*On the other hand, by the Lipschitz gradient condition of $f$ [29] we have:*

$$\begin{aligned} f(\mathbf{x}_{k+1}) \leq & f(\mathbf{x}_k) + \langle \nabla f(\mathbf{x}_k), \mathbf{x}_{k+1} - \mathbf{x}_k \rangle + \frac{L}{2}\|\mathbf{x}_{k+1} - \mathbf{x}_k\|^2 \\ \leq & f(\mathbf{x}) - \frac{\mu}{2}\|\mathbf{x} - \mathbf{x}_k\|^2 + \langle \nabla f(\mathbf{x}_k), \mathbf{x}_{k+1} - \mathbf{x} \rangle + \frac{L}{2}\|\mathbf{x}_{k+1} - \mathbf{x}_k\|^2. \end{aligned}$$

*So*

$$\begin{aligned} F(\mathbf{x}_{k+1}) - F(\mathbf{x}) \leq & \left\langle \mathbf{A}^T \overline{\boldsymbol{\lambda}}_{k+1}, \mathbf{x} - \mathbf{x}_{k+1} \right\rangle + \eta_k \left\langle \mathbf{x}_{k+1} - \mathbf{x}_k + \sigma_k, \mathbf{x} - \mathbf{x}_{k+1} \right\rangle \\ & + \frac{L}{2}\|\mathbf{x}_{k+1} - \mathbf{x}_k\|^2 - \frac{\mu}{2}\|\mathbf{x} - \mathbf{x}_k\|^2 + \epsilon_k \\ = & \left\langle \mathbf{A}^T \overline{\boldsymbol{\lambda}}_{k+1}, \mathbf{x} - \mathbf{x}_{k+1} \right\rangle + \eta_k \left\langle \mathbf{x}_{k+1} - \mathbf{x}_k, \mathbf{x} - \mathbf{x}_k \right\rangle + \eta_k \left\langle \sigma_k, \mathbf{x} - \mathbf{x}_{k+1} \right\rangle \\ & - \frac{\mu}{2}\|\mathbf{x} - \mathbf{x}_k\|^2 - \left( \frac{L}{2} + \frac{\beta\|\mathbf{A}^T \mathbf{A}\|_2}{\alpha_k} \right) \|\mathbf{x}_{k+1} - \mathbf{x}_k\|^2 + \epsilon_k. \end{aligned}$$

*Letting $\mathbf{x} = \mathbf{x}_k$ and $\mathbf{x} = \mathbf{x}^*$, we have*

$$\begin{aligned} F(\mathbf{x}_{k+1}) - F(\mathbf{x}_k) \leq & \left\langle \mathbf{A}^T \overline{\boldsymbol{\lambda}}_{k+1}, \mathbf{x}_k - \mathbf{x}_{k+1} \right\rangle - \left( \frac{L}{2} + \frac{\beta\|\mathbf{A}^T \mathbf{A}\|_2}{\alpha_k} \right) \|\mathbf{x}_{k+1} - \mathbf{x}_k\|^2 \\ & + \eta_k \left\langle \sigma_k, \mathbf{x}_k - \mathbf{x}_{k+1} \right\rangle + \epsilon_k, \end{aligned}$$

*and*

$$\begin{aligned} F(\mathbf{x}_{k+1}) - F(\mathbf{x}^*) \leq & \left\langle \mathbf{A}^T \overline{\boldsymbol{\lambda}}_{k+1}, \mathbf{x}^* - \mathbf{x}_{k+1} \right\rangle + \eta_k \left\langle \mathbf{x}_{k+1} - \mathbf{x}_k, \mathbf{x}^* - \mathbf{x}_k \right\rangle \\ & + \eta_k \left\langle \sigma_k, \mathbf{x}^* - \mathbf{x}_{k+1} \right\rangle - \frac{\mu}{2}\|\mathbf{x}^* - \mathbf{x}_k\|^2 - \left( \frac{L}{2} + \frac{\beta\|\mathbf{A}^T \mathbf{A}\|_2}{\alpha_k} \right) \|\mathbf{x}_{k+1} - \mathbf{x}_k\|^2 + \epsilon_k. \end{aligned}$$

*Multiplying the first inequality by $1 - \theta_k$ and the second by $\theta_k$ and adding them, we have:*

$$F(\mathbf{x}_{k+1}) - (1 - \theta_k)F(\mathbf{x}_k) - \theta_k F(\mathbf{x}^*)$$
$$\leq \langle \overline{\boldsymbol{\lambda}}_{k+1}, \theta_k \mathbf{A}\mathbf{x}^* + (1 - \theta_k)\mathbf{A}\mathbf{x}_k - \mathbf{A}\mathbf{x}_{k+1} \rangle + \eta_k \theta_k \langle \mathbf{x}_{k+1} - \mathbf{x}_k, \mathbf{x}^* - \mathbf{x}_k \rangle$$
$$+ \eta_k \langle \sigma_k, \theta_k \mathbf{x}^* + (1 - \theta_k)\mathbf{x}_k - \mathbf{x}_{k+1} \rangle - \frac{\mu \theta_k}{2} \|\mathbf{x}^* - \mathbf{x}_k\|^2$$
$$- \left( \frac{L}{2} + \frac{\beta \|\mathbf{A}^T \mathbf{A}\|_2}{\alpha_k} \right) \|\mathbf{x}_{k+1} - \mathbf{x}_k\|^2 + \epsilon_k$$
$$\leq \langle \overline{\boldsymbol{\lambda}}_{k+1}, \theta_k \mathbf{A}\mathbf{x}^* + (1 - \theta_k)\mathbf{A}\mathbf{x}_k - \mathbf{A}\mathbf{x}_{k+1} \rangle + \eta_k \theta_k \langle \mathbf{x}_{k+1} - \mathbf{x}_k, \mathbf{x}^* - \mathbf{x}_k \rangle$$
$$- \frac{\mu \theta_k}{2} \|\mathbf{x}^* - \mathbf{x}_k\|^2 - \left( \frac{L}{2} + \frac{\beta \|\mathbf{A}^T \mathbf{A}\|_2}{\alpha_k} \right) \|\mathbf{x}_{k+1} - \mathbf{x}_k\|^2$$
$$+ \sqrt{2\epsilon_k \eta_k}(\|\mathbf{x}^* - \mathbf{x}_{k+1}\| + \|\mathbf{x}^* - \mathbf{x}_k\|) + \epsilon_k.$$

where we use

$$\eta_k \langle \sigma_k, \theta_k \mathbf{x}^* + (1 - \theta_k)\mathbf{x}_k - \mathbf{x}_{k+1} \rangle = \eta_k \langle \sigma_k, \mathbf{x}^* - \mathbf{x}_{k+1} \rangle + \eta_k (\theta_k - 1) \langle \sigma_k, \mathbf{x}^* - \mathbf{x}_k \rangle$$
$$\leq \eta_k \|\sigma_k\| \|\mathbf{x}^* - \mathbf{x}_{k+1}\| + \eta_k \|\sigma_k\| \|\mathbf{x}^* - \mathbf{x}_k\| \leq \sqrt{2\epsilon_k \eta_k}(\|\mathbf{x}^* - \mathbf{x}_{k+1}\| + \|\mathbf{x}^* - \mathbf{x}_k\|),$$

and the first inequality is due to $0 \leq \theta_k \leq 1$ and the second is from $\|\sigma_k\| \leq \sqrt{\frac{2\epsilon_k}{\eta_k}}$. So we have

$$F(\mathbf{x}_{k+1}) - F(\mathbf{x}^*) + \langle \boldsymbol{\lambda}^*, \mathbf{A}\mathbf{x}_{k+1} - \mathbf{b} \rangle - (1 - \theta_k)\left( F(\mathbf{x}_k) - F(\mathbf{x}^*) + \langle \boldsymbol{\lambda}^*, \mathbf{A}\mathbf{x}_k - \mathbf{b} \rangle \right)$$
$$= F(\mathbf{x}_{k+1}) - (1 - \theta_k)F(\mathbf{x}_k) - \theta_k F(\mathbf{x}^*) + \langle \boldsymbol{\lambda}^*, \mathbf{A}\mathbf{x}_{k+1} - (1 - \theta_k)\mathbf{A}\mathbf{x}_k - \theta_k \mathbf{b} \rangle$$
$$\leq \langle \overline{\boldsymbol{\lambda}}_{k+1} - \boldsymbol{\lambda}^*, \theta_k \mathbf{A}\mathbf{x}^* + (1 - \theta_k)\mathbf{A}\mathbf{x}_k - \mathbf{A}\mathbf{x}_{k+1} \rangle + \eta_k \theta_k \langle \mathbf{x}_{k+1} - \mathbf{x}_k, \mathbf{x}^* - \mathbf{x}_k \rangle$$
$$- \frac{\mu \theta_k}{2} \|\mathbf{x}^* - \mathbf{x}_k\|^2 - \left( \frac{L}{2} + \frac{\beta \|\mathbf{A}^T \mathbf{A}\|_2}{\alpha_k} \right) \|\mathbf{x}_{k+1} - \mathbf{x}_k\|^2$$
$$+ \sqrt{2\epsilon_k \eta_k}(\|\mathbf{x}_{k+1} - \mathbf{x}^*\| + \|\mathbf{x}_k - \mathbf{x}^*\|) + \epsilon_k$$
$$= \langle \overline{\boldsymbol{\lambda}}_{k+1} - \boldsymbol{\lambda}^*, \theta_k \mathbf{A}\mathbf{x}^* + (1 - \theta_k)\mathbf{A}\mathbf{x}_k - \mathbf{A}\mathbf{x}_{k+1} \rangle + \frac{\eta_k \theta_k}{2}\left( \|\mathbf{x}_k - \mathbf{x}^*\|^2 - \|\mathbf{x}_{k+1} - \mathbf{x}^*\|^2 \right)$$
$$- \frac{\mu \theta_k}{2} \|\mathbf{x}^* - \mathbf{x}_k\|^2 - \left( \frac{L}{2} + \frac{\beta \|\mathbf{A}^T \mathbf{A}\|_2}{\alpha_k} - \frac{\eta_k \theta_k}{2} \right) \|\mathbf{x}_{k+1} - \mathbf{x}_k\|^2$$
$$+ \sqrt{2\epsilon_k \eta_k}(\|\mathbf{x}_{k+1} - \mathbf{x}^*\| + \|\mathbf{x}_k - \mathbf{x}^*\|) + \epsilon_k$$
$$= \frac{\alpha_k}{\beta} \langle \overline{\boldsymbol{\lambda}}_{k+1} - \boldsymbol{\lambda}^*, \boldsymbol{\lambda}_k - \boldsymbol{\lambda}_{k+1} \rangle + \frac{\eta_k \theta_k}{2}\left( \|\mathbf{x}_k - \mathbf{x}^*\|^2 - \|\mathbf{x}_{k+1} - \mathbf{x}^*\|^2 \right)$$
$$- \frac{\mu \theta_k}{2} \|\mathbf{x}^* - \mathbf{x}_k\|^2 - \left( \frac{L}{2} + \frac{\beta \|\mathbf{A}^T \mathbf{A}\|_2}{\alpha_k} - \frac{\eta_k \theta_k}{2} \right) \|\mathbf{x}_{k+1} - \mathbf{x}_k\|^2$$
$$+ \sqrt{2\epsilon_k \eta_k}(\|\mathbf{x}_{k+1} - \mathbf{x}^*\| + \|\mathbf{x}_k - \mathbf{x}^*\|) + \epsilon_k$$
$$= \frac{\alpha_k}{2\beta}\left( \|\boldsymbol{\lambda}_k - \boldsymbol{\lambda}^*\|^2 - \|\boldsymbol{\lambda}_{k+1} - \boldsymbol{\lambda}^*\|^2 - \|\overline{\boldsymbol{\lambda}}_{k+1} - \boldsymbol{\lambda}_k\|^2 + \|\boldsymbol{\lambda}_{k+1} - \overline{\boldsymbol{\lambda}}_{k+1}\|^2 \right)$$
$$+ \frac{\eta_k \theta_k}{2}\left( \|\mathbf{x}_k - \mathbf{x}^*\|^2 - \|\mathbf{x}_{k+1} - \mathbf{x}^*\|^2 \right)$$
$$- \frac{\mu \theta_k}{2} \|\mathbf{x}^* - \mathbf{x}_k\|^2 - \left( \frac{L}{2} + \frac{\beta \|\mathbf{A}^T \mathbf{A}\|_2}{\alpha_k} - \frac{\eta_k \theta_k}{2} \right) \|\mathbf{x}_{k+1} - \mathbf{x}_k\|^2$$
$$+ \sqrt{2\epsilon_k \eta_k}(\|\mathbf{x}_{k+1} - \mathbf{x}^*\| + \|\mathbf{x}_k - \mathbf{x}^*\|) + \epsilon_k$$

$$\leq \frac{\alpha_k}{2\beta} \left( \|\boldsymbol{\lambda}_k - \boldsymbol{\lambda}^*\|^2 - \|\boldsymbol{\lambda}_{k+1} - \boldsymbol{\lambda}^*\|^2 - \|\overline{\boldsymbol{\lambda}}_{k+1} - \boldsymbol{\lambda}_k\|^2 \right)$$
$$+ \frac{\eta_k \theta_k}{2} \left( \|\mathbf{x}_k - \mathbf{x}^*\|^2 - \|\mathbf{x}_{k+1} - \mathbf{x}^*\|^2 \right)$$
$$- \frac{\mu \theta_k}{2} \|\mathbf{x}^* - \mathbf{x}_k\|^2 - \left( \frac{L}{2} + \frac{\beta \|\mathbf{A}^T \mathbf{A}\|_2}{2\alpha_k} - \frac{\eta_k \theta_k}{2} \right) \|\mathbf{x}_{k+1} - \mathbf{x}_k\|^2$$
$$+ \sqrt{2\epsilon_k \eta_k}(\|\mathbf{x}_{k+1} - \mathbf{x}^*\| + \|\mathbf{x}_k - \mathbf{x}^*\|) + \epsilon_k,$$

where we define $\boldsymbol{\lambda}_k = \frac{\beta}{\alpha_{k-1}}(\mathbf{A}\mathbf{x}_k - \mathbf{b})$ and use

$$\begin{aligned}
\boldsymbol{\lambda}_{k+1} - \boldsymbol{\lambda}_k &= \frac{\beta}{\alpha_k}(\mathbf{A}\mathbf{x}_{k+1} - \mathbf{b}) - \frac{\beta}{\alpha_{k-1}}(\mathbf{A}\mathbf{x}_k - \mathbf{b}) \\
&= \frac{\beta}{\alpha_k}(\mathbf{A}\mathbf{x}_{k+1} - \mathbf{b}) - \frac{\beta(1-\theta_k)}{\alpha_k}(\mathbf{A}\mathbf{x}_k - \mathbf{b}) \\
&= \frac{\beta}{\alpha_k}\left(\mathbf{A}\mathbf{x}_{k+1} - (1-\theta_k)\mathbf{A}\mathbf{x}_k - \theta_k \mathbf{b}\right)
\end{aligned} \tag{6}$$

in the third equality and

$$\frac{\alpha_k}{2\beta}\|\boldsymbol{\lambda}_{k+1} - \overline{\boldsymbol{\lambda}}_{k+1}\|^2 = \frac{\alpha_k}{2\beta}\left\|\frac{\beta}{\alpha_k}\mathbf{A}(\mathbf{x}_{k+1} - \mathbf{x}_k)\right\|^2 \leq \frac{\beta\|\mathbf{A}^T\mathbf{A}\|_2}{2\alpha_k}\|\mathbf{x}_{k+1} - \mathbf{x}_k\|^2 \tag{7}$$

in the second inequality. From $\frac{1-\theta_k}{\alpha_k} = \frac{1}{\alpha_{k-1}}$, we have $\theta_k = 1 - \frac{\alpha_k}{\alpha_{k-1}}$ and

$$\eta_k \theta_k = L - \frac{L\alpha_k}{\alpha_{k-1}} + \beta\|\mathbf{A}^T\mathbf{A}\|_2\left(\frac{1}{\alpha_k} - \frac{1}{\alpha_{k-1}}\right). \tag{8}$$

Thus we have $\frac{L}{2} + \frac{\beta\|\mathbf{A}^T\mathbf{A}\|_2}{2\alpha_k} - \frac{\eta_k \theta_k}{2} \geq 0$ and

$$F(\mathbf{x}_{k+1}) - F(\mathbf{x}^*) + \langle \boldsymbol{\lambda}^*, \mathbf{A}\mathbf{x}_{k+1} - \mathbf{b}\rangle - (1-\theta_k)\left(F(\mathbf{x}_k) - F(\mathbf{x}^*) + \langle \boldsymbol{\lambda}^*, \mathbf{A}\mathbf{x}_k - \mathbf{b}\rangle\right)$$
$$\leq \frac{\alpha_k}{2\beta}\left(\|\boldsymbol{\lambda}_k - \boldsymbol{\lambda}^*\|^2 - \|\boldsymbol{\lambda}_{k+1} - \boldsymbol{\lambda}^*\|^2\right) + \frac{(\eta_k - \mu)\theta_k}{2}\|\mathbf{x}_k - \mathbf{x}^*\|^2 - \frac{\eta_k \theta_k}{2}\|\mathbf{x}_{k+1} - \mathbf{x}^*\|^2$$
$$+ \sqrt{2\epsilon_k \eta_k}(\|\mathbf{x}_{k+1} - \mathbf{x}^*\| + \|\mathbf{x}_k - \mathbf{x}^*\|) + \epsilon_k.$$

Dividing both sides by $\alpha_k$ and using $\frac{1-\theta_k}{\alpha_k} = \frac{1}{\alpha_{k-1}}$ and (4) we have

$$\frac{1}{\alpha_k}\left(F(\mathbf{x}_{k+1}) - F(\mathbf{x}^*) + \langle \boldsymbol{\lambda}^*, \mathbf{A}\mathbf{x}_{k+1} - \mathbf{b}\rangle\right) - \frac{1}{\alpha_{k-1}}\left(F(\mathbf{x}_k) - F(\mathbf{x}^*) + \langle \boldsymbol{\lambda}^*, \mathbf{A}\mathbf{x}_k - \mathbf{b}\rangle\right)$$
$$\leq \frac{1}{2\beta}\left(\|\boldsymbol{\lambda}_k - \boldsymbol{\lambda}^*\|^2 - \|\boldsymbol{\lambda}_{k+1} - \boldsymbol{\lambda}^*\|^2\right) + \frac{(\eta_k - \mu)\theta_k}{2\alpha_k}\|\mathbf{x}_k - \mathbf{x}^*\|^2 - \frac{\eta_k \theta_k}{2\alpha_k}\|\mathbf{x}_{k+1} - \mathbf{x}^*\|^2$$
$$+ \frac{\sqrt{2\epsilon_k \eta_k}}{\alpha_k}(\|\mathbf{x}_{k+1} - \mathbf{x}^*\| + \|\mathbf{x}_k - \mathbf{x}^*\|) + \frac{\epsilon_k}{\alpha_k}$$
$$\leq \frac{1}{2\beta}\left(\|\boldsymbol{\lambda}_k - \boldsymbol{\lambda}^*\|^2 - \|\boldsymbol{\lambda}_{k+1} - \boldsymbol{\lambda}^*\|^2\right) + \frac{\eta_{k-1}\theta_{k-1}}{2\alpha_{k-1}}\|\mathbf{x}_k - \mathbf{x}^*\|^2 - \frac{\eta_k \theta_k}{2\alpha_k}\|\mathbf{x}_{k+1} - \mathbf{x}^*\|^2$$
$$+ \frac{\sqrt{2\epsilon_k \eta_k}}{\alpha_k}(\|\mathbf{x}_{k+1} - \mathbf{x}^*\| + \|\mathbf{x}_k - \mathbf{x}^*\|) + \frac{\epsilon_k}{\alpha_k}.$$

Summing over $k = 0, 1, \cdots, K$, we have

$$\frac{1}{\alpha_K}\left(F(\mathbf{x}_{K+1}) - F(\mathbf{x}^*) + \langle \boldsymbol{\lambda}^*, \mathbf{A}\mathbf{x}_{K+1} - \mathbf{b}\rangle + \frac{\eta_K \theta_K}{2}\|\mathbf{x}_{K+1} - \mathbf{x}^*\|^2\right)$$
$$+ \frac{1}{2\beta}\|\boldsymbol{\lambda}_{K+1} - \boldsymbol{\lambda}^*\|^2$$
$$\leq \frac{1}{2\beta}\|\boldsymbol{\lambda}_0 - \boldsymbol{\lambda}^*\|^2 + \sum_{k=0}^{K} \frac{\sqrt{2\epsilon_k \eta_k}}{\alpha_k}(\|\mathbf{x}_{k+1} - \mathbf{x}^*\| + \|\mathbf{x}_k - \mathbf{x}^*\|) + \sum_{k=0}^{K}\frac{\epsilon_k}{\alpha_k}$$
$$= \frac{1}{2\beta}\|\boldsymbol{\lambda}_0 - \boldsymbol{\lambda}^*\|^2 + \frac{\sqrt{2\epsilon_0 \eta_0}}{\alpha_0}\|\mathbf{x}_0 - \mathbf{x}^*\| + \sum_{k=0}^{K}\frac{2\sqrt{\epsilon_k}}{\sqrt{\theta_k}\alpha_k}\frac{\sqrt{\eta_k \theta_k}}{\sqrt{2\alpha_k}}\|\mathbf{x}_{k+1} - \mathbf{x}^*\|$$
$$+ \sum_{k=1}^{K}\frac{2\sqrt{\epsilon_k \eta_k \alpha_{k-1}}}{\alpha_k \sqrt{\eta_{k-1}\theta_{k-1}}}\frac{\sqrt{\eta_{k-1}\theta_{k-1}}}{\sqrt{2\alpha_{k-1}}}\|\mathbf{x}_k - \mathbf{x}^*\| + \sum_{k=0}^{K}\frac{\epsilon_k}{\alpha_k},$$

where we use $\frac{1}{\alpha_{-1}} = 0$, $\eta_{-1} = L$.

On the other hand, since $(\mathbf{x}^*, \boldsymbol{\lambda}^*)$ is the saddle point, we have $\mathbf{x}^* = \mathrm{argmin}_{\mathbf{x}} F(\mathbf{x}) + \langle \boldsymbol{\lambda}^*, \mathbf{A}\mathbf{x} - \mathbf{b}\rangle$. So

$$F(\mathbf{x}^*) + \langle \boldsymbol{\lambda}^*, \mathbf{A}\mathbf{x}^* - \mathbf{b}\rangle \leq F(\mathbf{x}_{k+1}) + \langle \boldsymbol{\lambda}^*, \mathbf{A}\mathbf{x}_{k+1} - \mathbf{b}\rangle.$$

Thus we have $0 \leq F(\mathbf{x}_{k+1}) - F(\mathbf{x}^*) + \langle \boldsymbol{\lambda}^*, \mathbf{A}\mathbf{x}_{k+1} - \mathbf{b}\rangle$. Let
$u_{k+1} = \frac{1}{\alpha_k}\left(F(\mathbf{x}_{k+1}) - F(\mathbf{x}^*) + \langle \boldsymbol{\lambda}^*, \mathbf{A}\mathbf{x}_{k+1} - \mathbf{b}\rangle + \frac{\eta_k \theta_k}{2}\|\mathbf{x}_{k+1} - \mathbf{x}^*\|^2\right) + \frac{1}{2\beta}\|\boldsymbol{\lambda}_{k+1} - \boldsymbol{\lambda}^*\|^2$, then

$$u_{K+1} \leq \frac{1}{2\beta}\|\boldsymbol{\lambda}_0 - \boldsymbol{\lambda}^*\|^2 + \frac{\sqrt{2\epsilon_0\eta_0}}{\alpha_0}\|\mathbf{x}_0 - \mathbf{x}^*\| + \sum_{k=0}^{K}\frac{2\sqrt{\epsilon_k}}{\sqrt{\theta_k}\alpha_k}\sqrt{u_{k+1}} + \sum_{k=1}^{K+1}\frac{2\sqrt{\epsilon_k \eta_k \alpha_{k-1}}}{\alpha_k \sqrt{\eta_{k-1}\theta_{k-1}}}\sqrt{u_k} + \sum_{k=0}^{K}\frac{\epsilon_k}{\alpha_k}$$
$$= S_{K+1} + \sum_{k=0}^{K}\left(\frac{2\sqrt{\epsilon_{k+1}\eta_{k+1}\alpha_k}}{\alpha_{k+1}\sqrt{\eta_k \theta_k}} + \frac{2\sqrt{\epsilon_k}}{\sqrt{\theta_k}\alpha_k}\right)\sqrt{u_{k+1}},$$

where we let $S_{k+1} = \frac{1}{2\beta}\|\boldsymbol{\lambda}_0 - \boldsymbol{\lambda}^*\|^2 + \frac{\sqrt{2\epsilon_0\eta_0}}{\alpha_0}\|\mathbf{x}_0 - \mathbf{x}^*\| + \sum_{k=0}^{K}\frac{\epsilon_k}{\alpha_k}$. Since $u_0 = \frac{1}{2\beta}\|\boldsymbol{\lambda}_0 - \boldsymbol{\lambda}^*\|^2 \leq S_0$, then from the following Lemma:

**Lemma 1** [31] Assume $\{S_k\}$ is increasing and $\{v_k\}, \{\alpha_i\}$ are nonnegative, $v_0^2 \leq S_0$. If $v_k^2 \leq S_k + \sum_{i=1}^{k}\alpha_i v_i$, then $v_k \leq \frac{1}{2}\sum_{i=1}^{k}\alpha_i + \sqrt{\left(\frac{1}{2}\sum_{i=1}^{k}\alpha_i\right)^2 + S_k}$.

We have

$$\sqrt{u_{K+1}} \leq \sqrt{S_{K+1}} + \sum_{k=0}^{K}\left(\frac{2\sqrt{\epsilon_{k+1}\eta_{k+1}\alpha_k}}{\alpha_{k+1}\sqrt{\eta_k \theta_k}} + \frac{2\sqrt{\epsilon_k}}{\sqrt{\theta_k}\alpha_k}\right),$$

and

$$u_{K+1} \leq 2S_{K+1} + 2\left(\sum_{k=0}^{K}\left(\frac{2\sqrt{\epsilon_{k+1}\eta_{k+1}\alpha_k}}{\alpha_{k+1}\sqrt{\eta_k \theta_k}} + \frac{2\sqrt{\epsilon_k}}{\sqrt{\theta_k}\alpha_k}\right)\right)^2$$
$$= \frac{1}{\beta}\|\boldsymbol{\lambda}_0 - \boldsymbol{\lambda}^*\|^2 + \frac{2\sqrt{2\epsilon_0\eta_0}}{\alpha_0}\|\mathbf{x}_0 - \mathbf{x}^*\| + 2\sum_{k=0}^{K}\frac{\epsilon_k}{\alpha_k} + 2\left(\sum_{k=0}^{K}\left(\frac{2\sqrt{\epsilon_{k+1}\eta_{k+1}\alpha_k}}{\alpha_{k+1}\sqrt{\eta_k \theta_k}} + \frac{2\sqrt{\epsilon_k}}{\sqrt{\theta_k}\alpha_k}\right)\right)^2 \equiv C.$$

So $F(\mathbf{x}_{K+1}) - F(\mathbf{x}^*) + \langle \boldsymbol{\lambda}^*, \mathbf{A}\mathbf{x}_{K+1} - \mathbf{b}\rangle \leq C\alpha_K$ and $\frac{1}{2\beta}\|\boldsymbol{\lambda}_{K+1} - \boldsymbol{\lambda}^*\|^2 \leq C$. Since $\left\|\frac{\beta}{\alpha_K}(\mathbf{A}\mathbf{x}_{K+1} - \mathbf{b})\right\| = \|\boldsymbol{\lambda}_{K+1}\| \leq \|\boldsymbol{\lambda}_{K+1} - \boldsymbol{\lambda}^*\| + \|\boldsymbol{\lambda}^*\| \leq \sqrt{2\beta C} + \|\boldsymbol{\lambda}^*\|$, then we have

$$\|\mathbf{A}\mathbf{x}_{K+1} - \mathbf{b}\| \leq \left(\frac{\sqrt{2\beta C} + \|\boldsymbol{\lambda}^*\|}{\beta}\right)\alpha_K.$$

Thus we have

$$F(\mathbf{x}_{K+1}) - F(\mathbf{x}^*) \leq C\alpha_K - \langle \boldsymbol{\lambda}^*, \mathbf{A}\mathbf{x}_{K+1} - \mathbf{b} \rangle \leq C\alpha_K + \|\boldsymbol{\lambda}^*\| \|\mathbf{A}\mathbf{x}_{K+1} - \mathbf{b}\|$$
$$\leq \left( C + \frac{\sqrt{2\beta C} + \|\boldsymbol{\lambda}^*\|}{\beta} \|\boldsymbol{\lambda}^*\| \right) \alpha_K,$$

and

$$F(\mathbf{x}_{K+1}) - F(\mathbf{x}^*) \geq -\|\boldsymbol{\lambda}^*\| \|\mathbf{A}\mathbf{x}_{K+1} - \mathbf{b}\| \geq -\|\boldsymbol{\lambda}^*\| \left( \frac{\sqrt{2\beta C} + \|\boldsymbol{\lambda}^*\|}{\beta} \right) \alpha_K,$$

which completes the proof.

Then we can specialize the choice of $\alpha_k$ at each iteration for the generally convex case and strongly convex case, respectively. Theorem 2 establishes the $O(1/\sqrt{K})$ and $O(1/K)$ convergence rates for these two cases, respectively.

**Theorem 2** *If the conditions in Theorem 1 hold and*

1. *if $\mu = 0$, $\alpha_k = \frac{1}{\sqrt{k+1}}$ and $\epsilon_k = \frac{1}{(k+1)^{5+\delta}}$, then assumptions (4) and (5) hold and we have*

   $$|F(\mathbf{x}_{K+1}) - F(\mathbf{x}^*)| \leq O\left(1/\sqrt{K}\right) \text{ and } \|\mathbf{A}\mathbf{x}_{K+1} - \mathbf{b}\| \leq O\left(1/\sqrt{K}\right).$$

2. *if $\mu > 0$, $\alpha_k = \theta_k = \frac{1}{k+1}$, $\beta \leq \frac{\mu}{\|\mathbf{A}^T\mathbf{A}\|_2}$ and $\epsilon_k = \frac{1}{(k+1)^{4+\delta}}$, then assumptions (4) and (5) hold and we have*

   $$|F(\mathbf{x}_{K+1}) - F(\mathbf{x}^*)| \leq O(1/K) \text{ and } \|\mathbf{A}\mathbf{x}_{K+1} - \mathbf{b}\| \leq O(1/K).$$

**Proof 2** *If $\mu = 0$ and $\alpha_k = \frac{1}{\sqrt{k+1}}$, then from (8) we have $\frac{\eta_k \theta_k}{\alpha_k} = L\left(\frac{1}{\alpha_k} - \frac{1}{\alpha_{k-1}}\right) + \frac{\beta \|\mathbf{A}^T\mathbf{A}\|_2}{\alpha_k}\left(\frac{1}{\alpha_k} - \frac{1}{\alpha_{k-1}}\right)$
$= \frac{L}{\sqrt{k+1}+\sqrt{k}} + \frac{\beta \|\mathbf{A}^T\mathbf{A}\|_2}{1+\sqrt{1-1/(k+1)}}$, which is decreasing as $k$ increases. Thus (4) holds. From $\frac{1-\theta_k}{\alpha_k} = \frac{1}{\alpha_{k-1}}$, we have $\theta_k = 1 - \frac{\sqrt{k}}{\sqrt{k+1}} = \frac{1}{\sqrt{k+1}(\sqrt{k+1}+\sqrt{k})}$. So*

$$\sum_{k=0}^{\infty} \frac{\epsilon_k}{\alpha_k} = \sum_{k=0}^{\infty} \frac{1}{(k+1)^{4.5+\delta}} < \infty,$$
$$\sum_{k=0}^{\infty} \left( \frac{2\sqrt{\epsilon_{k+1}\eta_{k+1}\alpha_k}}{\alpha_{k+1}\sqrt{\eta_k \theta_k}} + \frac{2\sqrt{\epsilon_k}}{\sqrt{\theta_k \alpha_k}} \right) = O\left( \sum_{k=0}^{\infty} \frac{\sqrt{\epsilon_k}}{\sqrt{\theta_k \alpha_k}} \right) = O\left( \sum_{k=0}^{\infty} \frac{1}{(k+1)^{1+\delta/2}} \right) < \infty.$$

*Thus (5) holds and we have $|F(\mathbf{x}_{K+1}) - F(\mathbf{x}^*)| \leq O\left(1/\sqrt{K}\right)$ and $\|\mathbf{A}\mathbf{x}_{K+1} - \mathbf{b}\| \leq O\left(1/\sqrt{K}\right)$.*

*If $\mu > 0$ and $\alpha_k = \theta_k$, then (4) reduces to $\eta_k - \mu \leq \eta_{k-1}$. Consider $\eta_k - \eta_{k-1} = \beta \|\mathbf{A}^T\mathbf{A}\|_2 \left(\frac{1}{\alpha_k} - \frac{1}{\alpha_{k-1}}\right) = \beta \|\mathbf{A}^T\mathbf{A}\|_2$, where we use $\frac{1}{\alpha_k} - \frac{1}{\alpha_{k-1}} = \frac{\theta_k}{\alpha_k} = 1$. So if $\beta \leq \frac{\mu}{\|\mathbf{A}^T\mathbf{A}\|_2}$, then (4) holds. On the other hand, from $\frac{1}{\alpha_k} - \frac{1}{\alpha_{k-1}} = 1$ and $\frac{1}{\alpha_{-1}} = 0$ we have $\alpha_k = \frac{1}{k+1}$.*

$$\sum_{k=0}^{\infty} \frac{\epsilon_k}{\alpha_k} = \sum_{k=1}^{\infty} \frac{1}{(k+1)^{3+\delta}} < \infty,$$
$$\sum_{k=0}^{\infty} \left( \frac{2\sqrt{\epsilon_{k+1}\eta_{k+1}\alpha_k}}{\alpha_{k+1}\sqrt{\eta_k \theta_k}} + \frac{2\sqrt{\epsilon_k}}{\sqrt{\theta_k \alpha_k}} \right) = O\left( \sum_{k=0}^{\infty} \frac{\sqrt{\epsilon_k}}{\sqrt{\theta_k \alpha_k}} \right) = O\left( \sum_{k=0}^{\infty} \frac{1}{(k+1)^{1+\delta/2}} \right) < \infty.$$

*Thus (5) holds and we have $|F(\mathbf{x}_{K+1}) - F(\mathbf{x}^*)| \leq O(1/K)$ and $\|\mathbf{A}\mathbf{x}_{K+1} - \mathbf{b}\| \leq O(1/K)$.*

## 2.2 APG as the Solver

In this section, we consider to use APG for problem (2). At each iteration, it consists of two steps:

$$\begin{aligned}
\mathbf{y}_k &= \mathbf{x}_k + \frac{(\eta_k \theta_k - \mu)(1 - \theta_{k-1})}{(\eta_k - \mu)\theta_{k-1}}(\mathbf{x}_k - \mathbf{x}_{k-1}), \\
\mathbf{x}_{k+1} &= \mathrm{prox}_{h,\epsilon_k}^{1/\eta_k}\left(\mathbf{y}_k - \frac{\nabla f(\mathbf{y}_k) + \frac{\beta}{\alpha_k}\mathbf{A}^T(\mathbf{A}\mathbf{y}_k - \mathbf{b})}{\eta_k}\right).
\end{aligned} \qquad (9)$$

Similar to Theorem 1, we first establish the following general theorem.

**Theorem 3** *Assume that $f(\mathbf{x})$ is strongly convex with modulus $\mu \geq 0$. Let $\alpha_k \geq 0$ be a decreasing sequence with $\frac{1}{\alpha_{-1}} = 0$. Define $\theta_k$ as $\frac{1 - \theta_k}{\alpha_k} = \frac{1}{\alpha_{k-1}}$ and $\eta_k = L + \frac{\beta \|\mathbf{A}^T\mathbf{A}\|_2}{\alpha_k}$. Assume*

$$\frac{\eta_k (\theta_k)^2 - \mu \theta_k}{2\alpha_k} \leq \frac{\eta_{k-1}(\theta_{k-1})^2}{2\alpha_{k-1}}, \quad \theta_k \geq \frac{\mu}{\eta_k}, \qquad (10)$$

$$\sum_{k=0}^{\infty} \frac{\epsilon_k}{\alpha_k} < \infty, \quad \sum_{k=0}^{\infty} \frac{2\sqrt{\epsilon_k}}{\sqrt{\alpha_k}} < \infty. \qquad (11)$$

*Then for APG based quadratic penalty method (9), we have*

$$|F(\mathbf{x}_{K+1}) - F(\mathbf{x}^*)| \leq O(\alpha_K), \quad \|\mathbf{A}\mathbf{x}_{K+1} - \mathbf{b}\| \leq O(\alpha_K),$$

*where $\mathbf{x}^*$ is an optimal solution to problem (1).*

**Proof 3** *Let $\overline{\boldsymbol{\lambda}}_{k+1} = \frac{\beta}{\alpha_k}(\mathbf{A}\mathbf{y}_k - \mathbf{b})$ and $\boldsymbol{\lambda}_{k+1} = \frac{\beta}{\alpha_k}(\mathbf{A}\mathbf{x}_{k+1} - \mathbf{b})$. Similar to Theorem 1, we can have*

$$h(\mathbf{x}) - h(\mathbf{x}_{k+1}) \geq -\epsilon_k - \langle \nabla f(\mathbf{y}_k) + \mathbf{A}^T \overline{\boldsymbol{\lambda}}_{k+1} + \eta_k(\mathbf{x}_{k+1} - \mathbf{y}_k + \sigma_k), \mathbf{x} - \mathbf{x}_{k+1} \rangle,$$

$$f(\mathbf{x}_{k+1}) \leq f(\mathbf{x}) - \frac{\mu}{2}\|\mathbf{x} - \mathbf{y}_k\|^2 + \langle \nabla f(\mathbf{y}_k), \mathbf{x}_{k+1} - \mathbf{x} \rangle + \frac{L}{2}\|\mathbf{x}_{k+1} - \mathbf{y}_k\|^2.$$

*So*

$$\begin{aligned}
F(\mathbf{x}_{k+1}) - F(\mathbf{x}) \leq\ & \langle \mathbf{A}^T \overline{\boldsymbol{\lambda}}_{k+1}, \mathbf{x} - \mathbf{x}_{k+1} \rangle + \eta_k \langle \mathbf{x}_{k+1} - \mathbf{y}_k, \mathbf{x} - \mathbf{y}_k \rangle + \eta_k \langle \sigma_k, \mathbf{x} - \mathbf{x}_{k+1} \rangle \\
& - \frac{\mu}{2}\|\mathbf{x} - \mathbf{y}_k\|^2 - \left(\frac{L}{2} + \frac{\beta \|\mathbf{A}^T\mathbf{A}\|_2}{\alpha_k}\right)\|\mathbf{x}_{k+1} - \mathbf{y}_k\|^2 + \epsilon_k,
\end{aligned}$$

*Letting $\mathbf{x} = \mathbf{x}_k$ and $\mathbf{x} = \mathbf{x}^*$, we have*

$$\begin{aligned}
F(\mathbf{x}_{k+1}) - F(\mathbf{x}_k) \leq\ & \langle \mathbf{A}^T \overline{\boldsymbol{\lambda}}_{k+1}, \mathbf{x}_k - \mathbf{x}_{k+1} \rangle + \eta_k \langle \mathbf{x}_{k+1} - \mathbf{y}_k, \mathbf{x}_k - \mathbf{y}_k \rangle + \eta_k \langle \sigma_k, \mathbf{x}_k - \mathbf{x}_{k+1} \rangle \\
& - \left(\frac{L}{2} + \frac{\beta \|\mathbf{A}^T\mathbf{A}\|_2}{\alpha_k}\right)\|\mathbf{x}_{k+1} - \mathbf{y}_k\|^2 + \epsilon_k, \\
F(\mathbf{x}_{k+1}) - F(\mathbf{x}^*) \leq\ & \langle \mathbf{A}^T \overline{\boldsymbol{\lambda}}_{k+1}, \mathbf{x}^* - \mathbf{x}_{k+1} \rangle + \eta_k \langle \mathbf{x}_{k+1} - \mathbf{y}_k, \mathbf{x}^* - \mathbf{y}_k \rangle + \eta_k \langle \sigma_k, \mathbf{x}^* - \mathbf{x}_{k+1} \rangle \\
& - \frac{\mu}{2}\|\mathbf{x}^* - \mathbf{y}_k\|^2 - \left(\frac{L}{2} + \frac{\beta \|\mathbf{A}^T\mathbf{A}\|_2}{\alpha_k}\right)\|\mathbf{x}_{k+1} - \mathbf{y}_k\|^2 + \epsilon_k,
\end{aligned}$$

*Similar to Theorem 1, we have*

$$
\begin{aligned}
&F(\mathbf{x}_{k+1}) - F(\mathbf{x}^*) + \langle \boldsymbol{\lambda}^*, \mathbf{A}\mathbf{x}_{k+1} - \mathbf{b} \rangle - (1 - \theta_k)\left(F(\mathbf{x}_k) - F(\mathbf{x}^*) + \langle \boldsymbol{\lambda}^*, \mathbf{A}\mathbf{x}_k - \mathbf{b} \rangle\right) \\
&\leq \langle \overline{\boldsymbol{\lambda}}_{k+1} - \boldsymbol{\lambda}^*, \theta_k \mathbf{A}\mathbf{x}^* + (1 - \theta_k)\mathbf{A}\mathbf{x}_k - \mathbf{A}\mathbf{x}_{k+1} \rangle \\
&\quad + \eta_k \langle \mathbf{x}_{k+1} - \mathbf{y}_k, \theta_k \mathbf{x}^* + (1 - \theta_k)\mathbf{x}_k - \mathbf{y}_k \rangle \\
&\quad - \frac{\mu \theta_k}{2}\|\mathbf{x}^* - \mathbf{y}_k\|^2 - \left(\frac{L}{2} + \frac{\beta\|\mathbf{A}^T\mathbf{A}\|_2}{\alpha_k}\right)\|\mathbf{x}_{k+1} - \mathbf{y}_k\|^2 \\
&\quad + \eta_k \langle \sigma_k, \theta_k \mathbf{x}^* + (1 - \theta_k)\mathbf{x}_k - \mathbf{x}_{k+1} \rangle + \epsilon_k \\
&= \frac{\alpha_k}{\beta} \langle \overline{\boldsymbol{\lambda}}_{k+1} - \boldsymbol{\lambda}^*, \boldsymbol{\lambda}_k - \boldsymbol{\lambda}_{k+1} \rangle \\
&\quad + \frac{\eta_k}{2} \left(\|\theta_k \mathbf{x}^* + (1 - \theta_k)\mathbf{x}_k - \mathbf{y}_k\|^2 - \|\theta_k \mathbf{x}^* + (1 - \theta_k)\mathbf{x}_k - \mathbf{x}_{k+1}\|^2\right) \\
&\quad - \frac{\mu \theta_k}{2}\|\mathbf{x}^* - \mathbf{y}_k\|^2 - \frac{\beta\|\mathbf{A}^T\mathbf{A}\|_2}{2\alpha_k}\|\mathbf{y}_k - \mathbf{x}_{k+1}\|^2 \\
&\quad + \eta_k \langle \sigma_k, \theta_k \mathbf{x}^* + (1 - \theta_k)\mathbf{x}_k - \mathbf{x}_{k+1} \rangle + \epsilon_k \\
&= \frac{\alpha_k}{2\beta} \left(\|\boldsymbol{\lambda}_k - \boldsymbol{\lambda}^*\|^2 - \|\boldsymbol{\lambda}_{k+1} - \boldsymbol{\lambda}^*\|^2 - \|\overline{\boldsymbol{\lambda}}_{k+1} - \boldsymbol{\lambda}_k\|^2 + \|\boldsymbol{\lambda}_{k+1} - \overline{\boldsymbol{\lambda}}_{k+1}\|^2\right) \\
&\quad + \frac{\eta_k (\theta_k)^2}{2} \left(\left\|\mathbf{x}^* + \frac{1 - \theta_k}{\theta_k}\mathbf{x}_k - \frac{1}{\theta_k}\mathbf{y}_k\right\|^2 - \left\|\mathbf{x}^* + \frac{1 - \theta_k}{\theta_k}\mathbf{x}_k - \frac{1}{\theta_k}\mathbf{x}_{k+1}\right\|^2\right) \\
&\quad - \frac{\mu \theta_k}{2}\|\mathbf{x}^* - \mathbf{y}_k\|^2 - \frac{\beta\|\mathbf{A}^T\mathbf{A}\|_2}{2\alpha_k}\|\mathbf{y}_k - \mathbf{x}_{k+1}\|^2 \\
&\quad + \eta_k \theta_k \left\langle \sigma_k, \mathbf{x}^* + \frac{1 - \theta_k}{\theta_k}\mathbf{x}_k - \frac{1}{\theta_k}\mathbf{x}_{k+1} \right\rangle + \epsilon_k \\
&\leq \frac{\alpha_k}{2\beta} \left(\|\boldsymbol{\lambda}_k - \boldsymbol{\lambda}^*\|^2 - \|\boldsymbol{\lambda}_{k+1} - \boldsymbol{\lambda}^*\|^2 - \|\overline{\boldsymbol{\lambda}}_{k+1} - \boldsymbol{\lambda}_k\|^2\right) - \frac{\mu \theta_k}{2}\|\mathbf{x}^* - \mathbf{y}_k\|^2 \\
&\quad + \frac{\eta_k (\theta_k)^2}{2} \left(\left\|\mathbf{x}^* + \frac{1 - \theta_k}{\theta_k}\mathbf{x}_k - \frac{1}{\theta_k}\mathbf{y}_k\right\|^2 - \left\|\mathbf{x}^* + \frac{1 - \theta_k}{\theta_k}\mathbf{x}_k - \frac{1}{\theta_k}\mathbf{x}_{k+1}\right\|^2\right), \\
&\quad + \eta_k \theta_k \left\langle \sigma_k, \mathbf{x}^* + \frac{1 - \theta_k}{\theta_k}\mathbf{x}_k - \frac{1}{\theta_k}\mathbf{x}_{k+1} \right\rangle + \epsilon_k.
\end{aligned}
$$

*where we use (6) in the first equality and (7) in the last inequality but replacing $\mathbf{x}_k$ with $\mathbf{y}_k$ in (7). Consider*

$$
\begin{aligned}
&\frac{\eta_k (\theta_k)^2}{2} \left\|\mathbf{x}^* + \frac{1 - \theta_k}{\theta_k}\mathbf{x}_k - \frac{1}{\theta_k}\mathbf{y}_k\right\|^2 \\
&= \frac{\eta_k (\theta_k)^2}{2} \left\|\frac{\mu}{\eta_k \theta_k}(\mathbf{x}^* - \mathbf{y}_k) + \left(1 - \frac{\mu}{\eta_k \theta_k}\right)\left(\mathbf{x}^* + \frac{\eta_k(1 - \theta_k)}{\eta_k \theta_k - \mu}\mathbf{x}_k - \frac{\eta_k - \mu}{\eta_k \theta_k - \mu}\mathbf{y}_k\right)\right\|^2 \\
&\leq \frac{\mu \theta_k}{2}\|\mathbf{x}^* - \mathbf{y}_k\|^2 + \frac{\theta_k(\eta_k \theta_k - \mu)}{2}\left\|\mathbf{x}^* + \frac{\eta_k(1 - \theta_k)}{\eta_k \theta_k - \mu}\mathbf{x}_k - \frac{\eta_k - \mu}{\eta_k \theta_k - \mu}\mathbf{y}_k\right\|^2,
\end{aligned}
$$

*where we let $\theta_k \geq \frac{\mu}{\eta_k}$. Let $\mathbf{d}_{k+1} = \frac{1}{\theta_k}\mathbf{x}_{k+1} - \frac{1 - \theta_k}{\theta_k}\mathbf{x}_k$, $\mathbf{d}_k = \frac{\eta_k - \mu}{\eta_k \theta_k - \mu}\mathbf{y}_k - \frac{\eta_k(1 - \theta_k)}{\eta_k \theta_k - \mu}\mathbf{x}_k = \frac{1}{\theta_{k-1}}\mathbf{x}_k - $*

$\frac{1-\theta_{k-1}}{\theta_{k-1}}\mathbf{x}_{k-1}$, we have $\mathbf{y}_k = \mathbf{x}_k + \frac{(\eta_k\theta_k - \mu)(1-\theta_{k-1})}{(\eta_k - \mu)\theta_{k-1}}(\mathbf{x}_k - \mathbf{x}_{k-1})$ and

$$F(\mathbf{x}_{k+1}) - F(\mathbf{x}^*) + \langle \boldsymbol{\lambda}^*, \mathbf{A}\mathbf{x}_{k+1} - \mathbf{b}\rangle - (1-\theta_k)\left(F(\mathbf{x}_k) - F(\mathbf{x}^*) + \langle \boldsymbol{\lambda}^*, \mathbf{A}\mathbf{x}_k - \mathbf{b}\rangle\right)$$
$$\leq \frac{\alpha_k}{2\beta}\left(\|\boldsymbol{\lambda}_k - \boldsymbol{\lambda}^*\|^2 - \|\boldsymbol{\lambda}_{k+1} - \boldsymbol{\lambda}^*\|^2\right)$$
$$+ \frac{\eta_k(\theta_k)^2 - \mu\theta_k}{2}\|\mathbf{d}_k - \mathbf{x}^*\|^2 - \frac{\eta_k(\theta_k)^2}{2}\|\mathbf{d}_{k+1} - \mathbf{x}^*\|^2$$
$$- \eta_k\theta_k\langle \sigma_k, \mathbf{d}_{k+1} - \mathbf{x}^*\rangle + \epsilon_k$$
$$\leq \frac{\alpha_k}{2\beta}\left(\|\boldsymbol{\lambda}_k - \boldsymbol{\lambda}^*\|^2 - \|\boldsymbol{\lambda}_{k+1} - \boldsymbol{\lambda}^*\|^2\right)$$
$$+ \frac{\eta_k(\theta_k)^2 - \mu\theta_k}{2}\|\mathbf{d}_k - \mathbf{x}^*\|^2 - \frac{\eta_k(\theta_k)^2}{2}\|\mathbf{d}_{k+1} - \mathbf{x}^*\|^2$$
$$+ \sqrt{2\eta_k\epsilon_k}\theta_k\|\mathbf{d}_{k+1} - \mathbf{x}^*\| + \epsilon_k.$$

Then following the same proof in Theorem 1 and using $\frac{\eta_k(\theta_k)^2 - \mu\theta_k}{2\alpha_k} \leq \frac{\eta_{k-1}(\theta_{k-1})^2}{2\alpha_{k-1}}$, we have

$$\frac{1}{\alpha_k}\left(F(\mathbf{x}_{k+1}) - F(\mathbf{x}^*) + \langle\boldsymbol{\lambda}^*, \mathbf{A}\mathbf{x}_{k+1} - \mathbf{b}\rangle + \frac{\eta_k(\theta_k)^2}{2}\|\mathbf{d}_{k+1} - \mathbf{x}^*\|^2\right)$$
$$- \frac{1}{\alpha_{k-1}}\left(F(\mathbf{x}_k) - F(\mathbf{x}^*) + \langle\boldsymbol{\lambda}^*, \mathbf{A}\mathbf{x}_k - \mathbf{b}\rangle + \frac{\eta_{k-1}(\theta_{k-1})^2}{2}\|\mathbf{d}_k - \mathbf{x}^*\|^2\right)$$
$$+ \frac{1}{2\beta}\|\boldsymbol{\lambda}_{k+1} - \boldsymbol{\lambda}^*\|^2 - \frac{1}{2\beta}\|\boldsymbol{\lambda}_k - \boldsymbol{\lambda}^*\|^2$$
$$\leq \frac{\theta_k\sqrt{2\epsilon_k\eta_k}}{\alpha_k}\|\mathbf{d}_{k+1} - \mathbf{x}^*\| + \frac{\epsilon_k}{\alpha_k}.$$

Summing over $k = 0, \cdots, K$, we have

$$\frac{1}{\alpha_k}\left(F(\mathbf{x}_{K+1}) - F(\mathbf{x}^*) + \langle\boldsymbol{\lambda}^*, \mathbf{A}\mathbf{x}_{K+1} - \mathbf{b}\rangle + \frac{\eta_K(\theta_K)^2}{2}\|\mathbf{d}_{K+1} - \mathbf{x}^*\|^2\right)$$
$$+ \frac{1}{2\beta}\|\boldsymbol{\lambda}_{k+1} - \boldsymbol{\lambda}^*\|^2$$
$$\leq \frac{1}{2\beta}\|\boldsymbol{\lambda}_0 - \boldsymbol{\lambda}^*\|^2 + \sum_{k=0}^K \frac{2\sqrt{\epsilon_k}}{\sqrt{\alpha_k}}\frac{\sqrt{\eta_k}\theta_k}{\sqrt{2\alpha_k}}\|\mathbf{d}_{k+1} - \mathbf{x}^*\| + \sum_{k=0}^K \frac{\epsilon_k}{\alpha_k},$$

Similar to Theorem 1, we have

$$\frac{1}{\alpha_k}\left(F(\mathbf{x}_{K+1}) - F(\mathbf{x}^*) + \langle\boldsymbol{\lambda}^*, \mathbf{A}\mathbf{x}_{K+1} - \mathbf{b}\rangle + \frac{\eta_K(\theta_K)^2}{2}\|\mathbf{d}_{K+1} - \mathbf{x}^*\|^2\right)$$
$$+ \frac{1}{2\beta}\|\boldsymbol{\lambda}_{k+1} - \boldsymbol{\lambda}^*\|^2$$
$$\leq \frac{1}{\beta}\|\boldsymbol{\lambda}_0 - \boldsymbol{\lambda}^*\|^2 + 2\sum_{k=0}^K \frac{\epsilon_k}{\alpha_k} + 2\left(\sum_{k=0}^K \frac{2\sqrt{\epsilon_k}}{\sqrt{\alpha_k}}\right)^2.$$

Following the same proof as Theorem 1, we can have the conclusion.

Then we can specialize the value of $\alpha_k$ for the generally convex case and the strongly convex case and establish the convergence rates in Theorem 4.

**Theorem 4** *If the conditions in Theorem 3 hold and*

1. if $\mu = 0$, $\alpha_k = \theta_k = \frac{1}{k+1}$ and $\epsilon_k = \frac{1}{(k+1)^{3+\delta}}$, then assumptions (10) and (11) hold and we have

    $|F(\mathbf{x}_{K+1}) - F(\mathbf{x}^*)| \leq O(1/K)$ and $\|\mathbf{A}\mathbf{x}_{K+1} - \mathbf{b}\| \leq O(1/K)$.

2. if $\mu > 0$, $\frac{1-\theta_k}{(\theta_k)^2} = \frac{1}{(\theta_{k-1})^2}$, $\alpha_k = (\theta_k)^2$, $\theta_0 = 1$, $\frac{\mu^2}{4L\|\mathbf{A}^T\mathbf{A}\|} \leq \beta \leq \frac{\mu}{\|\mathbf{A}^T\mathbf{A}\|_2}$ and $\epsilon_k = \frac{1}{(k+1)^{4+\delta}}$, then assumptions (10) and (11) hold and we have

    $|F(\mathbf{x}_{K+1}) - F(\mathbf{x}^*)| \leq O(1/K^2)$ and $\|\mathbf{A}\mathbf{x}_{K+1} - \mathbf{b}\| \leq O(1/K^2)$.

**Proof 4** If $\mu = 0$ and $\alpha_k = \theta_k$, then (10) reduces to $\eta_k \theta_k \leq \eta_{k-1}\theta_{k-1}$ and $\theta_k \geq 0$, which is true due to $0 \leq \theta_k < \theta_{k-1}$ and the definition of $\eta_k$. From $\frac{1-\theta_k}{\alpha_k} = \frac{1}{\alpha_{k-1}}$ and $\frac{1}{\alpha_{-1}} = 0$ we have $\alpha_k = \frac{1}{k+1}$.

$$\sum_{k=0}^{\infty} \frac{\epsilon_k}{\alpha_k} = \sum_{k=0}^{\infty} \frac{1}{(k+1)^{2+\delta}} < \infty,$$

$$\sum_{k=0}^{\infty} \frac{2\sqrt{\epsilon_k}}{\sqrt{\alpha_k}} = \sum_{k=0}^{\infty} \frac{1}{(k+1)^{1+\delta/2}} < \infty.$$

So we have $\|\mathbf{A}\mathbf{x}_{K+1} - \mathbf{b}\| \leq O(1/K)$ and $|F(\mathbf{x}_{K+1}) - F(\mathbf{x}^*)| \leq O(1/K)$.

If $\mu > 0$ and $\alpha_k = (\theta_k)^2$, then (10) reduces to $\eta_k - \mu/\theta_k \leq \eta_{k-1}$ and $L\theta_k + \frac{\beta\|\mathbf{A}^T\mathbf{A}\|_2}{\theta_k} \geq \mu$. Consider

$\eta_k - \mu/\theta_k - \eta_{k-1} = L + \frac{\beta\|\mathbf{A}^T\mathbf{A}\|_2}{\alpha_k} - \mu/\theta_k - \left(L + \frac{\beta\|\mathbf{A}^T\mathbf{A}\|_2}{\alpha_{k-1}}\right)$

$= \beta\|\mathbf{A}^T\mathbf{A}\|_2 \left(\frac{1}{\alpha_k} - \frac{1}{\alpha_{k-1}}\right) - \mu/\theta_k = \frac{\beta\|\mathbf{A}^T\mathbf{A}\|_2 - \mu}{\theta_k}$, where we use $\frac{1}{\alpha_k} - \frac{1}{\alpha_{k-1}} = \frac{\theta_k}{\alpha_k} = \frac{1}{\theta_k}$. So if $\beta \leq \frac{\mu}{\|\mathbf{A}^T\mathbf{A}\|_2}$, then $\eta_k - \mu/\theta_k \leq \eta_{k-1}$.

Then we consider $L\theta_k + \frac{\beta\|\mathbf{A}^T\mathbf{A}\|_2}{\theta_k} \geq \mu$. It holds if $\beta \geq \frac{\theta_k\mu - L(\theta_k)^2}{\|\mathbf{A}^T\mathbf{A}\|_2}$. Since $\theta\mu - L\theta^2 \leq \frac{\mu^2}{4L}, \forall \theta$, we only need $\beta \geq \frac{\mu^2}{4L\|\mathbf{A}^T\mathbf{A}\|_2}$. So we finally get the condition $\frac{\mu^2}{4L\|\mathbf{A}^T\mathbf{A}\|_2} \leq \beta \leq \frac{\mu}{\|\mathbf{A}^T\mathbf{A}\|_2}$. From $\theta_0 = 1$ and $\frac{1-\theta_k}{(\theta_k)^2} = \frac{1}{(\theta_{k-1})^2}$, we can easily have $\theta_k \leq \frac{2}{k+2}$ and $\alpha_k \leq \frac{4}{(k+2)^2}$.

$$\sum_{k=0}^{\infty} \frac{\epsilon_k}{\alpha_k} = O\left(\sum_{k=0}^{\infty} \frac{1}{(k+1)^{2+\delta}}\right) < \infty,$$

$$\sum_{k=0}^{\infty} \frac{2\sqrt{\epsilon_k}}{\sqrt{\alpha_k}} = O\left(\sum_{k=0}^{\infty} \frac{1}{(k+2)^{1+\delta/2}}\right) < \infty.$$

Thus we have $\|\mathbf{A}\mathbf{x}_{K+1} - \mathbf{b}\| \leq O(1/K^2)$ and $|F(\mathbf{x}_{K+1}) - F(\mathbf{x}^*)| \leq O(1/K^2)$.

## 3 Counterexamples for Using the Usual Solvers

In this section, we give counterexamples to show that using the usual PG or APG for computing the inexact solution to (2) may not converge to the solution to (1). This justifies the special choice of the sequences of penalty parameters shown in Sections 2.1 and 2.2 for the strongly convex case.

### 3.1 A Counterexample for the Usual PG Solver

In (1), we choose $h(\mathbf{x}) = f(\mathbf{x}) = \frac{1}{2}\|\mathbf{x}\|^2$, $\mathbf{A} = (\mathbf{D}\ \mathbf{0}) \in \mathbb{R}^{n \times 2n}$, and $\mathbf{b}$ is any $n \times 1$ vector, where $\mathbf{D} = \text{diag}(d_1, \cdots, d_n)$ is an $n \times n$ diagonal matrix with $|d_1| \geq \cdots \geq |d_n|$ and $d_1 \neq 0$. More specifically, our exemplar problem is:

$$\min_{\mathbf{x}} \frac{1}{2}\|\mathbf{x}\|^2 + \frac{1}{2}\|\mathbf{x}\|^2, \quad s.t. \quad (\mathbf{D}\ \mathbf{0})\mathbf{x} = \mathbf{b}, \tag{12}$$

and (3) can be explicitly written as

$$\mathbf{x}_{k+1} = \frac{1}{2\alpha_k \beta^{-1} d_1^{-2} + 1} \begin{pmatrix} \mathbf{I} - d_1^{-2}\mathbf{D}^2 & \mathbf{0} \\ \mathbf{0} & \mathbf{I} \end{pmatrix} \mathbf{x}_k + \frac{1}{2\alpha_k \beta^{-1} + d_1^2} \begin{pmatrix} \mathbf{D}\mathbf{b} \\ \mathbf{0} \end{pmatrix}.$$

The latter $n$ entries of $\mathbf{x}_{k+1}$ can be written down as:

$$\mathbf{x}_{k+1,i} = \frac{1}{2\alpha_k \beta^{-1} d_1^{-2} + 1}\mathbf{x}_{k,i} = \frac{1}{\prod\limits_{t=0}^{k}(2\alpha_t \beta^{-1} d_1^{-2} + 1)}\mathbf{x}_{0,i}, \quad i = n+1,\cdots,2n.$$

It is easy to see that the latter $n$ entries of the optimal solution to (12) are all zeros. if $\sum\limits_{t=0}^{+\infty}\alpha_t < +\infty$, then

$$\prod_{t=0}^{+\infty}(2\alpha_t\beta^{-1}d_1^{-2}+1) \leq \prod_{t=0}^{+\infty}\exp(2\alpha_t\beta^{-1}d_1^{-2}) = \exp\left(2\beta^{-1}d_1^{-2}\sum_{t=0}^{+\infty}\alpha_t\right) < +\infty,$$

and $\mathbf{x}_{k,i} > 0, \forall k, i = n+1,\cdots,2n$. So iterate (3) cannot converge to the optimal solution when $\alpha_k$ decreases faster than $O(1/k^{1+\delta})$, where $\delta$ is a small constant.

## 3.2 A Counterexample for the Usual APG Solver

We use counterexample (12) and set $\theta_k = \frac{2}{k+1}$, $\mu = 0$ in (9) for simplicity, which leads to $\mathbf{y}_k = \mathbf{x}_k + \frac{k-2}{k+1}(\mathbf{x}_k - \mathbf{x}_{k-1})$. The coefficient $\frac{k-2}{k+1}$ is adopted in [29] and is widely used in practice. Through the example, we have that

$$\mathbf{x}_{k+1,i} = \frac{1}{2\alpha_k\beta^{-1}d_1^{-2}+1}\mathbf{y}_{k,i} = (1-\tau_k)\left(\mathbf{x}_{k,i} + \frac{k-2}{k+1}(\mathbf{x}_{k,i} - \mathbf{x}_{k-1,i})\right), \forall i \geq n+1,$$

where we set $\tau_k = 1 - \frac{1}{2\alpha_k\beta^{-1}d_1^{-2}+1}$. Let $b_k = \mathbf{x}_{k,i} - \mathbf{x}_{k-1,i}$ and $\mathbf{x}_{1,i} = \mathbf{x}_{0,i} = \frac{1}{\delta} + 1$ with $0 < \delta < 0.7$, then $b_{k+1} = \frac{k-2}{k+1}(1-\tau_k)b_k - \tau_k\mathbf{x}_{k,i}$. We can use induction to prove that when $0 < \tau_k \leq \frac{1}{(k+1)^{2+\delta}\mathbf{x}_{0,i}}$, then $0 \geq b_k \geq -\frac{1}{(k+1)^{1+\delta}}$ and $\mathbf{x}_{k,i} > 1, \forall k$.

For $k = 2$, since $\mathbf{x}_{1,i} = \mathbf{x}_{0,i}$, so $b_1 = 0$, $b_2 = -\tau_1\mathbf{x}_{0,i} \leq 0$ and $b_2 = -\tau_1\mathbf{x}_{0,i} \geq -\frac{1}{2^{2+\delta}} \geq -\frac{1}{3^{1+\delta}}$ when $0 < \delta < 0.7$, which satisfies the condition.

Assume that $0 \geq b_k \geq -\frac{1}{(k+1)^{1+\delta}}$ holds for $k \leq t$, we have $\mathbf{x}_{k,i} - \mathbf{x}_{k-1,i} \geq -\frac{1}{(k+1)^{1+\delta}}, \forall k \leq t$. Summing with $k=1$ to $t$, we have $\mathbf{x}_{t,i} - \mathbf{x}_{0,i} \geq -\sum_{j=1}^{t}\frac{1}{(j+1)^{1+\delta}}$. So $\mathbf{x}_{t,i} \geq \mathbf{x}_{0,i} - \sum_{j=1}^{t}\frac{1}{(j+1)^{1+\delta}} > \mathbf{x}_{0,i} - \int_1^\infty \frac{1}{x^{1+\delta}}dx = \frac{1}{\delta} + 1 - \frac{1}{\delta} = 1$. We also have $\mathbf{x}_{t,i} \leq \mathbf{x}_{t-1,i} \leq \cdots \leq \mathbf{x}_{0,i}$.

Since $0 < \tau_t < 1$, $b_t \leq 0$ and $\mathbf{x}_{t,i} > 1$, we can easily have $b_{t+1} \leq 0$.

On the other hand, $b_{t+1} \geq -\frac{t-2}{t+1}(1-\tau_t)\frac{1}{(t+1)^{1+\delta}} - \tau_t\mathbf{x}_{0,i} \geq -\frac{t-2}{t+1}\frac{1}{(t+1)^{1+\delta}} - \tau_t\mathbf{x}_{0,i}$. So we only need $\frac{1}{(t+2)^{1+\delta}} - \frac{t-2}{t+1}\frac{1}{(t+1)^{1+\delta}} \geq \tau_t\mathbf{x}_{0,i} = \frac{1}{(t+1)^{2+\delta}}$. It is equivalent to $\left(1 - \frac{1}{t+2}\right)^{1+\delta} \geq \frac{t-1}{t+1}$. Since $\left(1 - \frac{1}{t+2}\right)^{1+\delta} \geq 1 - \frac{1+\delta}{t+2}, \forall \delta \geq 0$, we only need $1 - \frac{1+\delta}{t+2} \geq \frac{t-1}{t+1}$, which is satisfied when $\frac{t+3}{t+1} \geq \delta$.

Since $\tau_k = \frac{2\alpha_k\beta^{-1}d_1^{-2}}{2\alpha_k\beta^{-1}d_1^{-2}+1} \leq 2\alpha_k\beta^{-1}d_1^{-2}$, then $\alpha_k \leq \frac{\beta d_1^2}{2(k+1)^{2+\delta}\mathbf{x}_{0,i}}$ makes iterate (9) not converge to the optimal solution.

## 4 Decentralized Distributed Optimization

In this section, we consider the decentralized consensus optimization, where $n$ agents form a connected network and cooperatively solve problem:

$$\min \frac{1}{n}\sum_i (f_i(\mathbf{x}) + h_i(\mathbf{x})). \tag{13}$$

Each $f_i$ and $h_i$ are only available to agent $i$ and a pair of agents can exchange information if and only if they are connected directly. We consider the distributed algorithms using only local computation and communication, i.e., each agent $i$ maintains a local variable $\mathbf{x}(i)$, locally computes the gradient of $f_i$ and proximal operation of $h_i$ to update $\mathbf{x}(i)$, and then performs one communication with its neighbors.

We reformulate problem (13) as (1). Let $\mathbf{x} = \begin{pmatrix} \mathbf{x}(1)^T \\ \vdots \\ \mathbf{x}(n)^T \end{pmatrix}$, $f(\mathbf{x}) = \sum_i f_i(\mathbf{x}(i))$, $h(\mathbf{x}) = \sum_i h_i(\mathbf{x}(i))$, $\nabla f(\mathbf{x}) = \begin{pmatrix} \nabla f_1(\mathbf{x}(1))^T \\ \vdots \\ \nabla f_n(\mathbf{x}(n))^T \end{pmatrix}$. Assume that $f_i$ has $L_i$-Lipschitz continuous gradient, then $f$ has $(\max_i L_i)$-Lipschitz continuous gradient. We only consider the general convex case in this section. The analysis can be directly extended to the strongly convex case. Formulate the problem as:

$$\min f(\mathbf{x}) + h(\mathbf{x}), \quad s.t. \quad \mathbf{Ux} = 0, \tag{14}$$

where $\mathbf{U} = \left(\frac{\mathbf{I}-\mathbf{W}}{2}\right)^{1/2}$ and $\mathbf{W}$ is a symmetric and double stochastic matrix. $\mathbf{W}_{i,j} > 0$ if and only if agents $i$ and $j$ are neighbors or $i = j$, otherwise, $\mathbf{W}_{i,j} = 0$. Such $\mathbf{W}$ satisfies null$(\mathbf{I} - \mathbf{W}) = \{\mathbf{1}\}$ and $1 = \lambda_1(\mathbf{W}) \geq \lambda_2(\mathbf{W}) \geq \cdots \geq \lambda_n(\mathbf{W}) \geq -1$, where $\lambda_i(\mathbf{W})$ is the $i$-th largest eigenvalue. Thus $\|\mathbf{U}\|_2 = 1$ and constraint $\mathbf{Ux} = 0$ ensures that $\mathbf{x}$ is consensual, i.e., $\mathbf{x}(1) = \cdots = \mathbf{x}(n)$.

We use the penalty method to solve problem (14). For the PG based penalty method, at each iteration, (3) becomes

$$\begin{aligned}
\mathbf{x}_{k+1} &= \text{prox}_{h,\epsilon_k}^{1/\eta_k}\left(\mathbf{x}_k - \frac{\nabla f(\mathbf{x}_k) + \frac{\beta}{\alpha_k}\mathbf{U}^T\mathbf{U}\mathbf{x}_k}{\eta_k}\right) \\
&= \text{prox}_{h,\epsilon_k}^{1/\eta_k}\left(\hat{\mathbf{W}}^k\mathbf{x}_k - \frac{1}{L+\beta/\alpha_k}\nabla f(\mathbf{x}_k)\right),
\end{aligned} \tag{15}$$

where $\hat{\mathbf{W}}^k = \frac{(L+\beta/2\alpha_k)\mathbf{I} + (\beta/2\alpha_k)\mathbf{W}}{L+\beta/\alpha_k} \rightarrow \frac{\mathbf{I}+\mathbf{W}}{2} \equiv \hat{\mathbf{W}}$ as $\alpha_k \rightarrow 0$ and $\eta_k = L + \frac{\beta}{\alpha_k}$. Then $\hat{\mathbf{W}}$ is also symmetric and double stochastic with $1 = \lambda_1(\hat{\mathbf{W}}) \geq \lambda_2(\hat{\mathbf{W}}) \geq \cdots \geq \lambda_n(\hat{\mathbf{W}}) \geq 0$. For each agent $i$, (15) can be described as:

1. Receive $\mathbf{x}_k(j)$ from its neighborhoods.
2. Compute $\mathbf{x}_{k+1}(i) = \text{prox}_{h_i,\epsilon_k/n}^{1/\eta_k}\left(\sum_j \hat{\mathbf{W}}_{i,j}^k \mathbf{x}_k(j) - \frac{1}{L+\beta/\alpha_k}\nabla f_i(\mathbf{x}_k(i))\right)$.
3. Broadcast $\mathbf{x}_{k+1}(i)$ to its neighborhoods.

We can find that this strategy becomes the distributed gradient method. Specially, [13] studied the iterate of $\mathbf{x}_{k+1}(i) = \sum_j \widetilde{\mathbf{W}}_{i,j}\mathbf{x}_k(j) - \alpha \nabla f_i(\mathbf{x}_k(i))$ [11] with a constant step size, where $\widetilde{\mathbf{W}}$ respects the adjacent structure of the network and is symmetric and double stochastic with $1 = \lambda_1(\widetilde{\mathbf{W}}) > \lambda_2(\widetilde{\mathbf{W}}) \geq \lambda_3(\widetilde{\mathbf{W}}) \geq \cdots \geq \lambda_n(\widetilde{\mathbf{W}}) > -1$. In this case, the algorithm does not converge to the solution of problem (13) but a point in its neighborhood. This fact corresponds to the common sense that a constant penalty makes the penalty method not converge to the optimal solution of the constrained problem [1]. [19] proposed a distributed dual averaging algorithm, which has the special case of $\mathbf{x}_{k+1}(i) = \sum_j \widetilde{\mathbf{W}}_{i,j}\mathbf{x}_k(j) - \alpha_k \nabla f_i(\mathbf{x}_k(i))$ with a decreasing step size of $\alpha_k = \frac{1}{\sqrt{k}}$. [19] established the ergodic $O\left(\frac{\log K}{\sqrt{K}}\right)$ convergence rate in the sense of $\sum_i f_i\left(\frac{1}{K}\sum_{k=1}^K \mathbf{x}_k(t)\right) - \sum_i f_i(\mathbf{x}^*) \leq O\left(\frac{\log K}{\sqrt{K}}\right), \forall t$. [12] studied a slight different scheme of $\mathbf{u}_k(i) = \sum_j \widetilde{\mathbf{W}}_{i,j}\mathbf{x}_k(j)$, $\mathbf{x}_{k+1}(i) =$ $\text{prox}_{h_i}^{\alpha_k}[\mathbf{u}_k(i) - \alpha_k \nabla f_i(\mathbf{u}_k(i))]$ with $\alpha_k = \frac{1}{\sqrt{k}}$ and also established the $O\left(\frac{\log K}{\sqrt{K}}\right)$ convergence rate of $\min_{1\leq k\leq K} \sum_i f_i(\mathbf{x}_k(t)) - \sum_i f_i(\mathbf{x}^*) \leq O\left(\frac{\log K}{\sqrt{K}}\right)$. The technique in [19, 12] is to analyze the dual averaging method (gradient descent method) with the deviation of $\left\|\mathbf{x}(i) - \frac{1}{n}\sum_j \mathbf{x}(j)\right\|$ and control the deviation by the diminishing step size. As a comparison, we analyze it from the view of the penalty method and accordingly we can improve the convergence rate to $O\left(\frac{1}{\sqrt{K}}\right)$ under the same algorithm framework. Moreover, we have less requirement on the network topology, i.e., [13, 19] requires the spectral gap of $1 - \max\{|\lambda_2(\widehat{\mathbf{W}})|, |\lambda_n(\widehat{\mathbf{W}})|\} > 0$ while we have no such assumption on $\mathbf{W}$. Specially, we have the following theorem.

**Theorem 5** *Assume that $f_i + h_i$ is $l_i$-Lipschitz continuous, $\forall i$. Let $\alpha_k = \frac{1}{\sqrt{k+1}}$ and $\epsilon_k = \frac{1}{(k+1)^{5+\delta}}$, then for the PG based penalty method, we have*

$$\left[\sum_i f_i(\mathbf{x}_k(t)) + h_i(\mathbf{x}_k(t))\right] - \left[\sum_i f_i(\mathbf{x}^*) + h_i(\mathbf{x}^*)\right] \leq O\left(\frac{1}{\sqrt{k}}\right),$$

*where $\mathbf{x}^*$ is a minimum of problem (13) and $\mathbf{x}_k(t)$ is the local variable in any agent $t$.*

**Proof 5** *From Theorem 2, we have*

$$|f(\mathbf{x}_k) + h(\mathbf{x}_k) - f(\mathbf{x}^*) - h(\mathbf{x}^*)| \leq O(1/\sqrt{k}), \quad \|\mathbf{U}\mathbf{x}_k\| \leq O(1/\sqrt{k}).$$

*Rewrite $\mathbf{x}_k$ as $\mathbf{x}_k = [\mathbf{b}_{:,1}, \cdots, \mathbf{b}_{:,d}]$. Since $\mathbf{W}$ is a symmetric and double stochastic matrix, then*

$$O(1/k) \geq \|\mathbf{U}\mathbf{x}_k\|^2 = \sum_{i=1}^d \mathbf{b}_{:,i}^T \mathbf{U}^2 \mathbf{b}_{:,i} = \frac{1}{2}\sum_{i=1}^d \mathbf{b}_{:,i}^T(\mathbf{I}-\mathbf{W})\mathbf{b}_{:,i}$$

$$=\frac{1}{4}\sum_{i=1}^d\sum_{p=1}^n\sum_{q=1}^n \mathbf{W}_{p,q}(\mathbf{b}_{p,i}-\mathbf{b}_{q,i})^2 = \frac{1}{4}\sum_{p=1}^n\sum_{q=1}^n \mathbf{W}_{p,q}\|\mathbf{b}_{p,:}-\mathbf{b}_{q,:}\|^2$$

$$=\frac{1}{4}\sum_{p=1}^n\sum_{q=1}^n \mathbf{W}_{p,q}\|\mathbf{x}_k(p)-\mathbf{x}_k(q)\|^2 \geq \frac{w_{min}}{4}\sum_{p,q,\mathbf{W}_{p,q}>0}\|\mathbf{x}_k(p)-\mathbf{x}_k(q)\|^2,$$

*where we let $w_{\min} = \min_{p,q,\mathbf{W}_{p,q}>0} \mathbf{W}_{p,q}$. Then*

$$O(1/\sqrt{k}) \geq \sum_{p,q,\mathbf{W}_{p,q}>0}\|\mathbf{x}_k(p)-\mathbf{x}_k(q)\|.$$

*Moreover, let $p_1$ and $p_r$ be two nodes such that they do not connect directly. Since the network is connected, then there exists a connected path $p_1, \cdots, p_r$ and*

$$O(1/\sqrt{k}) \geq \|\mathbf{x}_k(p_1)-\mathbf{x}_k(p_2)\| + \cdots + \|\mathbf{x}_k(p_{r-1})-\mathbf{x}_k(p_r)\| \geq \|\mathbf{x}_k(p_1)-\mathbf{x}_k(p_r)\|.$$

*So for any $t$, we have*

$$\sum_i [f_i(\mathbf{x}_k(t)) + h_i(\mathbf{x}_k(t)) - f_i(\mathbf{x}^*) - h_i(\mathbf{x}^*)]$$

$$= \sum_i [f_i(\mathbf{x}_k(t)) + h_i(\mathbf{x}_k(t)) - f_i(\mathbf{x}_k(i)) - h_i(\mathbf{x}_k(i))]$$

$$+ \sum_i [f_i(\mathbf{x}_k(i)) + h_i(\mathbf{x}_k(i)) - f_i(\mathbf{x}^*) - h_i(\mathbf{x}^*)]$$

$$\leq \sum_i l_i \|\mathbf{x}_k(t) - \mathbf{x}_k(i)\| + O(1/\sqrt{k})$$

$$= O(1/\sqrt{k}),$$

*which completes the proof.*

Now we use the APG based penalty method for problem (14), then (9) becomes

$$\mathbf{y}_k = \mathbf{x}_k + \frac{\theta_k(1-\theta_{k-1})}{\theta_{k-1}}(\mathbf{x}_k - \mathbf{x}_{k-1}),$$
$$\mathbf{x}_{k+1} = \text{prox}_{h,\epsilon_k}^{1/\eta_k}\left(\hat{\mathbf{W}}_k\mathbf{y}_k - \frac{1}{L+1/\alpha_k}\nabla f(\mathbf{y}_k)\right), \qquad (16)$$

At each iteration, each agent $i$ performs the following operations:

1. Receive $\mathbf{y}_k(j)$ from its neighborhoods.

2. Compute $\mathbf{x}_{k+1}(i) = \text{prox}_{h_i,\epsilon_k/n}^{1/\eta_k}\left(\sum_j \hat{\mathbf{W}}_{i,j}^k \mathbf{y}_k(j) - \frac{1}{L+1/\alpha_k}\nabla f_i(\mathbf{y}_k(i))\right)$
   and $\mathbf{y}_{k+1}(i) = \mathbf{x}_{k+1}(i) + \frac{\theta_{k+1}(1-\theta_k)}{\theta_k}(\mathbf{x}_{k+1}(i) - \mathbf{x}_k(i))$.

3. Broadcast $\mathbf{y}_{k+1}(i)$ to its neighborhoods.

Similar to the relation between the PG based penalty method and the distributed gradient method, we also find that iterate (16) becomes the fast gradient method [18], which has the iterate of $\mathbf{x}_{k+1}(i) = \sum_j \widetilde{\mathbf{W}}_{i,j}\mathbf{y}_k(j) - \alpha_k \nabla f_i(\mathbf{y}_k(i)), \mathbf{y}_{k+1}(i) = \mathbf{x}_{k+1}(i) + \frac{k+1}{k+4}(\mathbf{x}_{k+1}(i) - \mathbf{x}_k(i))$, where the network topology matrix $\widetilde{\mathbf{W}}$ is symmetric and double stochastic with $\widetilde{\mathbf{W}} \succ 0$ and $\lambda_2(\widetilde{\mathbf{W}}) < 1$. [18] established the convergence rate of $O\left(\frac{\log K}{K}\right)$ with a decreasing step size of $\alpha_k = \frac{1}{k}$ by the tool of inexact Nesterov's gradient method [32], where the inexactness comes from $\left\|\mathbf{x}(i) - \frac{1}{n}\sum_j \mathbf{x}(j)\right\|$ and it is controlled by the decreasing step size. Due to the totally different analysis framework from the view of the penalty method, we can improve the convergence rate to $O\left(\frac{1}{K}\right)$ and do not require the network topology satisfying $\mathbf{W} \succ 0$ and $\lambda_2(\mathbf{W}) < 1$. Similar to Theorem 5, we have the following theorem for the APG based penalty method.

**Theorem 6** *Assume that $f_i + h_i$ is $l_i$ continuous, $\forall i$. Let $\alpha_k = \frac{1}{k+1}$ and $\epsilon_k = \frac{1}{(k+1)^{3+\delta}}$, then for the APG based penalty method, we have*

$$\left[\sum_i f_i(\mathbf{x}_k(t)) + h_i(\mathbf{x}_k(t))\right] - \left[\sum_i f_i(\mathbf{x}^*) + h_i(\mathbf{x}^*)\right] \leq O\left(\frac{1}{k}\right),$$

*where $\mathbf{x}^*$ is a minimum of problem (13) and $\mathbf{x}_k(t)$ is the local variable in any agent $t$.*

EXTRA [10] is another decentralized first order method with $O\left(\frac{1}{K}\right)$ convergence rate in the sense of $\frac{1}{K}\sum_{k=1}^K \|\mathbf{U}\mathbf{x}_k + \nabla f(\mathbf{x}_k)\|^2 \leq O\left(\frac{1}{K}\right)$ and $\frac{1}{K}\sum_{k=1}^K \|\mathbf{U}\mathbf{x}_k\|^2 \leq O\left(\frac{1}{K}\right)$. However, for the constraint, they only established $\min_{k\leq K}\|\mathbf{U}\mathbf{x}_k\| \leq o\left(\frac{1}{\sqrt{K}}\right)$. As a comparison, we can establish $\|\mathbf{U}\mathbf{x}_K\| \leq O\left(\frac{1}{\sqrt{K}}\right)$ for the PG based penalty method and $\|\mathbf{U}\mathbf{x}_K\| \leq O\left(\frac{1}{K}\right)$ for the APG based penalty method. EXTRA [10] also requires $\lambda_n(\mathbf{W}) > -1$.

### 4.1 Communication Efficient Distributed Optimization

In this section, we consider the case that communication is a major bottleneck in distributed optimization. In these days, CPUs can read and write the memory at over 10 GB per second whereas communication over TCP/IP is about 10 MB per second [25]. Thus such a consideration is necessary. In this section, we consider the case of $f = 0$ and the proximal mapping of $h$ has no closed form solution. Then the fast gradient method in [18] cannot handle such a nonsmooth problem and

**Algorithm 1** AdaptSmooth [33]

    **for** $t = 0, 1, \cdots, T_k - 1$ **do**
        Define $G_{i,k}^{\gamma_t}(\mathbf{u})$ as (18).
        Apply the accelerated gradient method [35] to minimize $G_{i,k}^{\gamma_t}$ for $O(\frac{1}{\sqrt{\gamma_t \sigma}})$ iterations such that
$G_{i,k}^{\gamma_t}(\mathbf{u}_{t+1}) - \min G_{i,k}^{\gamma_t} \leq \frac{G_{i,k}^{\gamma_t}(\mathbf{u}_t) - \min G_{i,k}^{\gamma_t}}{4}$.
        $\gamma_{t+1} = \gamma_t/2$.
    **end for**
    Output $\mathbf{x}_{k+1}(i) = \mathbf{u}_{T_k}$.

---

the dual averaging method [19] needs $O\left(\frac{\log 1/\varepsilon}{\varepsilon^2}\right)$ communications and computations to reach an $\varepsilon$ solution. We consider to use the APG based penalty method to improve the results to $O\left(\frac{1}{\varepsilon}\right)$ communications and $O\left(\frac{1}{\varepsilon^2}\right)$ computations, matching the best known result of the primal-dual based decentralized method in [25].

We consider the following iterate:

$$\begin{aligned} \mathbf{y}_k &= \mathbf{x}_k + \frac{\theta_k(1 - \theta_{k-1})}{\theta_{k-1}}(\mathbf{x}_k - \mathbf{x}_{k-1}), \\ \mathbf{x}_{k+1} &= \mathrm{prox}_{h, \epsilon_k}^{1/\eta_k}\left(\hat{\mathbf{W}} \mathbf{y}_k\right), \end{aligned} \quad (17)$$

where $\hat{\mathbf{W}} = \frac{\mathbf{I} + \mathbf{W}}{2}$. Theorem 6 tells us that if we can control the error as $\epsilon_k = \frac{1}{(k+1)^{3+\delta}}$, then iterate (17) only needs $O\left(\frac{1}{\varepsilon}\right)$ outer iterations, i.e., $O\left(\frac{1}{\varepsilon}\right)$ communications. So we only need to consider the inner iterations in (17) to control the computation cost.

Define $G_{i,k}(\mathbf{x}(i)) = h_i(\mathbf{x}(i)) + \frac{\beta}{2\alpha_k}\|\mathbf{x}(i) - \sum_j \hat{\mathbf{W}}_{i,j} \mathbf{y}_k(j)\|^2$. Then the goal of each agent $i$ at the $k$-th iterations is to find an $\epsilon_k/n$ minimizer of $G_{i,k}$. $G_{i,k}$ is $\frac{\beta}{\alpha_k}$ strongly convex but nonsmooth. Many algorithms can be used to minimize $G_{i,k}$ and we use the AdaptSmooth method in [33], which has the best convergence rate as far as we know and is described in Algorithm 1. In the algorithm we define

$$\begin{aligned} h_i^\gamma(\mathbf{x}(i)) &= \max_{\mathbf{z}} \left[\langle \mathbf{z}, \mathbf{x}(i)\rangle - h_i^*(\mathbf{z}) - \frac{\gamma}{2}\|\mathbf{z}\|^2\right], \\ G_{i,k}^\gamma(\mathbf{x}(i)) &= h_i^\gamma(\mathbf{x}(i)) + \frac{\beta}{2\alpha_k}\left\|\mathbf{x}(i) - \sum_j \hat{\mathbf{W}}_{i,j} \mathbf{y}_k(j)\right\|^2. \end{aligned} \quad (18)$$

where $h_i^*$ is the Fenchel dual of $h_i$. Then $G_{i,k}^\gamma$ is a smooth approximation of $G_{i,k}$. It is $\frac{\beta}{\alpha_k} \equiv \sigma$ strongly convex and has $\left(\frac{1}{\gamma} + \frac{\beta}{\alpha_k}\right)$-Lipschitz continuous gradient. Theorem 4.2 in [33] states that if each $h_i$ is $l_i$-Lipschitz continuous and $G_{i,k}(\mathbf{u}_0) - \min G_{i,k} \leq \triangle$, then AdaptSmooth with $\gamma_0 = \triangle/l_i^2$ and $T_k = \log_2(\triangle/\epsilon_k)$ produces $\mathbf{u}_{T_k}$ such that $G_{i,k}(\mathbf{u}_{T_k}) - \min G_{i,k} \leq \epsilon_k$ in $\frac{1}{\sqrt{\epsilon_k \sigma}}$ total computations of $\nabla G_{i,k}^\gamma$. From the settings in Theorem 6, we know $\frac{1}{\sqrt{\epsilon_k \sigma}} = 1/\sqrt{\frac{\beta(k+1)}{(k+1)^{3+\delta}}} = O\left((k+1)^{1+\delta/2}\right)$. So iterate (17) needs $\sum_{k=1}^K (k+1)^{1+\delta/2} = O\left(K^{2+\delta}\right)$ total computations of $\nabla G_{i,k}^\gamma$ to reach an $O(1/K)$ solution, i.e., $O\left(\frac{1}{\varepsilon^{2+\delta}}\right)$ computations for an $\varepsilon$ solution. In general, evaluating $\nabla G_{i,k}^\gamma$ may be more costly than obtaining the subgradient of $h_i$. However, when $G_{i,k}^\gamma$ is given in closed form, smoothing is efficient. For more details, please see [34].

## 5 Conclusion

In this paper, we study the quadratic penalty method with continuation and inexact update for linearly constrained convex problems, i.e., performing the proximal gradient descent only once and then increase the penalty. For such a variant of the quadratic penalty method, we give counterexamples to show that it may not converge to the solution to the original constrained problem. By carefully choosing the penalty parameters, for the PG based penalty method, we establish the convergence rates of $O\left(\frac{1}{\sqrt{K}}\right)$ and $O\left(\frac{1}{K}\right)$ when $f$ is generally convex and strongly convex, respectively; for APG we improve the convergence rates to $O\left(\frac{1}{K}\right)$ and $O\left(\frac{1}{K^2}\right)$.

When applied to the decentralized distributed optimization, the penalty methods become the widely used distributed gradient method and fast distributed gradient method. However, due to the different analysis framework, we can improve their $O\left(\frac{\log K}{\sqrt{K}}\right)$ and $O\left(\frac{\log K}{K}\right)$ convergence rates to $O\left(\frac{1}{\sqrt{K}}\right)$ and $O\left(\frac{1}{K}\right)$ under fewer assumptions on the network topology for general convex problems. We also give a communication efficient version of the fast distributed gradient method, which can find an $\varepsilon$ solution in $O\left(\frac{1}{\varepsilon}\right)$ communications and $O\left(\frac{1}{\varepsilon^{2+\delta}}\right)$ computations for the non-smooth problems, where $\delta$ is a small constant.